\documentclass{amsart}
\usepackage{amssymb}\usepackage{amsmath}\usepackage{a4}
\usepackage{color}

\def\Ireg{\mathcal{I}_{\operatorname{Reg}}}
\def\Iregc{\mathcal{I}_{\operatorname{Reg,c}}}
\def\dx{dx}\def\dy{dy}
\def\Id{\operatorname{Id}}
\def\Tr{\operatorname{Tr}}\def\iD{{\sqrt{-1}D}}
\newcommand{\nn}{\nonumber}
\newcommand{\al}{{\vec\alpha}}\newcommand{\D}{{\Delta_M}}
\def\DD{{\mathcal{D}}}
\newcommand{\vard}{y,r,\omega ,D_r,\lambda}
\newcommand{\var}{y,r,\omega,\tau,\lambda}

\newcommand{\gt}{\tilde g ^{ab}}
\newcommand{\gtr}{\tilde g ^{ab}_{,r}}
\newcommand{\gtrr}{\tilde g ^{ab}_{,rr}}

\newcommand{\tgr}{\tilde g _{ab,r}}
\newcommand{\tgrr}{\tilde g _{ab,rr}}
\newcommand{\tl}{\tau^2 + \Lambda^2}

\begin{document}
\newtheorem{theorem}{Theorem}[section]
\newtheorem{lemma}[theorem]{Lemma}
\newtheorem{remark}[theorem]{Remark}
\newtheorem{example}[theorem]{Example}
\def\qedbox{\hbox{$\rlap{$\sqcap$}\sqcup$}}
\makeatletter
  \renewcommand{\theequation}{%
   \thesection.\alph{equation}}
  \@addtoreset{equation}{section}
 \makeatother
\def\Op{\operatorname{Op}}
\title[Heat trace asymptotics]
{Heat trace asymptotics with \\singular weight functions}
\author{M. van den Berg, P. Gilkey, K. Kirsten, and R. Seeley}
\begin{address}{MvdB: Department of Mathematics, University of
Bristol, University Walk, Bristol,\newline\phantom{...a}BS8 1TW, U.K.}\end{address}
\begin{email}{M.vandenBerg@bris.ac.uk}\end{email}
\begin{address}{PG: Mathematics Department, University of Oregon, Eugene, OR 97403, USA}\end{address}
\begin{email}{gilkey@uoregon.edu}\end{email}
\begin{address}{KK: Department of Mathematics, Baylor University \\ Waco, TX 76798, USA}\end{address}
\begin{email}{Klaus\_Kirsten@baylor.edu}\end{email}
\begin{address}{RS: 35 Lakewood Rd, Newton, MA 02461 USA}\end{address}
\begin{email}{r-seeley@comcast.net}\end{email}
\begin{abstract} We study the weighted heat trace asymptotics of an operator of Laplace type with Dirichlet boundary conditions
where the weight function exhibits radial blowup. We give formulas for the first few terms in the expansion in
terms of geometrical data.
\end{abstract}
\maketitle

\section{Introduction}\label{sect-1}
\subsection{Motivation} The asymptotic analysis of the heat trace provides a natural link
between the spectrum of Laplace type operators $\mathcal{D}$ acting
on functions on a $m$ dimensional Riemannian manifold $M$ and the
underlying geometry of $M$. For small time $t$ it links the
distribution of the large energy part of the spectrum of
$\mathcal{D}$ to local geometric invariants of $M$ and its boundary
which show up in its asymptotic expansion. These invariants play an
important in many physical phenomena e.g. in quantum statistical
mechanics when taking the large volume limits or in the Casimir
effect \cite{K02}. Typically the coefficient of the leading
$t^{-m/2}$ term in the heat trace expansion for small t is
determined by the interior (volume) of $M$. In many situations a
detailed study of the boundary behaviour of the heat kernel
associated with $\partial_t+\mathcal{D}$ is desirable. One way of
obtaining this information is putting a weight in the evaluation of
the heat trace. In the setting of the heat content of $M$ this
corresponds to giving $M$ a non uniform specific heat. It is well
known that the diagonal element at $(x,x;t)$ of the Dirichlet heat
kernel associated to $e^{-t\mathcal{D}}$ vanishes like $r^2$ where
$r=r(x)$ is the distance to the boundary. This allows the weights to
diverge like $r^{-\alpha}$ , where $\operatorname{Re}(\alpha) < 3$.
We will show, using pseudo differential calculus, that a modified
asymptotic series still exists in this case. For example if $1 <
\alpha < 3$ the leading behaviour of the heat trace is
$t^{(1-m-\alpha)/2}$ with a coefficient determined by an integral
over the boundary of $M$.

\subsection{The heat equation}\label{sect-1.2}
We adopt the {\it Einstein convention} and sum over repeated
indices. Let $M$ be a compact smooth Riemannian manifold of dimension
$m$ and with smooth boundary $\partial M$. Let $V$ be a smooth vector bundle over $M$ and let $\DD$ be an operator of
Laplace type on the space $C^\infty(V)$ of smooth sections to $V$. This means that the leading symbol of
$\DD$ is given by the metric tensor or, equivalently, that we may express in any system of
local coordinates
$x=(x^1,...,x^m)$ and relative to any local frame for $V$ the operator $\DD$ in the form:
\begin{equation}\label{eqn-1.a}
\DD=-\left\{g^{\mu\nu}\operatorname{Id}\partial_{x_\mu}\partial_{x_\nu}
+A_1^\nu\partial_{x_\nu}+A_0\right\}\,.
\end{equation}
In Equation (\ref{eqn-1.a}), we let $1\le\mu,\nu\le m$, we let $A_1^\nu$
and
$A_0$ be smooth endomorphisms (matrices), and we let $g^{\mu\nu}$ be the inverse of the metric
$g_{\mu\nu}:=g(\partial_{x_\mu},\partial_{x_\nu})$. Note that the Riemannian measure $dx$ on $M$ is given by
$$
\dx:=gdx^1...dx^m\quad\text{where}\quad g:=\sqrt{\det(g_{\mu\nu})}\,.
$$
Thus, for example, the scalar Laplacian $\Delta_M:=\delta d$ is of Laplace type since
\begin{equation}\label{eqn-1.b}
\Delta_M=-\left(g^{\mu\nu}\partial_{x_\nu}\partial_{x_\mu}
+g^{-1}\partial_{x_\nu}\left\{gg^{\mu\nu}\right\}\partial_{x_\mu}\right)\,.
\end{equation}

We shall use the {\it Dirichlet realization} of the operator $\DD$.
For $t>0$ and for $\phi\in L^2(V)$, the heat equation
$$(\partial_t+\DD)u(x;t)=0,\quad
  u(\cdot;t)|_{\partial M}=0,\quad
  \lim_{t\downarrow 0}u(\cdot;t)=\phi(\cdot)\ \text{in}\ L^2(V),
$$
has a solution $u=e^{-t\DD}\phi$ which is smooth in $(x;t)$. The operator $e^{-t\DD}$  has a
kernel
$p_\DD(x,\tilde x;t)$ which is smooth in
$(x,\tilde x;t)$ such that
$$u(x;t)=\int_Mp_\DD(x,\tilde x;t)\phi(\tilde x)d\tilde x\,.$$
In the case of the scalar Laplacian $\Delta_M$, there is a complete orthonormal basis $\{\phi_i\}$ for
$L^2(M)$ where the $\phi_i\in C^\infty(M)$ satisfy $\phi_i|_{\partial M}=0$ and $\Delta_M\phi_i=\lambda_i\phi_i$.
The corresponding Dirichlet heat kernel $p_M:=p_{\Delta_M}$ is given in terms of the {\it spectral resolution}
$\{\phi_i,\lambda_i\}$ via
$$
p_M(x,\tilde x;t)=\sum_ie^{-t\lambda_i}\phi_i(x){{\bar\phi}}_i(\tilde x)\,.
$$

\subsection{Heat trace asymptotics in the smooth setting}\label{sect-1.3}
We use the geodesic flow defined by the unit inward
normal vector field to define a diffeomorphism for some $\varepsilon>0$ between the collar
${\mathcal{C}_{\varepsilon}}:=\partial M\times[0,\varepsilon]$ and a neighborhood of the boundary in $M$ which identifies
$\partial M\times\{0\}$ with $\partial M$;  the curves
$r\rightarrow (y_0,r)$ for $r\in[0,\varepsilon]$ are then unit speed geodesics perpendicular to the
boundary and $r$ is the geodesic distance to the boundary.
Let $F\in C^\infty(M)$ be an auxiliary weight
function which is used for localization. On $\mathcal{C}_\varepsilon$, expand $F$ in a Taylor series
$$F(y,r)\sim\sum_{i=0}^\infty F_i(y)r^i\quad\text{where}\quad F_i=\left.{\textstyle\frac1{i!}}(\partial_r)^iF\right|_{r=0}\,.$$
Henceforth we shall let $\Tr$ denote the fiber trace and $\Tr_{L^2}$ denote the global $L^2$ trace. We then have:
\begin{equation}\label{eqn-1.c}
\Tr_{L^2}(Fe^{-t\DD})=\int_MF(x)\Tr\{p_\DD(x,x;t)\}\dx\,.
\end{equation}
We note for future reference that on the diagonal, the heat kernel $p_{\mathbb{R}^m}(x,x;t)$ for $\mathbb{R}^m$ and the heat
kernel
$p_H(x,x;t)$ on the half space $H:=\{x:x_1>0\}$ of the scalar Laplacian are given by
\begin{equation}\label{eqn-1.d}
p_{\mathbb{R}^m}(x,x;t)=(4\pi t)^{-m/2}\quad\text{and}\quad
p_H(x,x;t)=(4\pi t)^{-m/2}(1-e^{-r^2/t})\,.
\end{equation}
Let $dy$ be the Riemannian measure on the boundary. To simplify future expressions, we set
$$\mathcal{I}\{F\}=\int_MF\dx\quad\text{and}\quad\mathcal{I}^{bd}\{F\}=\int_{\partial M}F\dy\,.$$
We will also use the notation $\mathcal{I}\{Fd\nu\}[U]$ when it is necessary to specify the
domain of integration $U$ and/or the measure $d\nu$.

\begin{theorem}\label{thm-1.1}
Let $M$ be a compact smooth Riemannian manifold. Let $\DD$ be the Dirichlet realization of an operator of Laplace type. Let $F\in
C^\infty(M)$.
\begin{enumerate}
\item There is a complete asymptotic
expansion as
$t\downarrow0$ of the form:
$$
\displaystyle\Tr_{L^2}(Fe^{-t\DD})\sim t^{-m/2}\sum_{n=0}^\infty
t^na_{n}(F,\DD)
+t^{-(m-1)/2}\sum_{\ell=0}^\infty t^{\ell/2}a_\ell^{bd}(F,\DD)\,.
$$
\item There are local invariants $a_n=a_n(x,\DD)$ defined on $M$ and there are local invariants
$a_{\ell,i}^{bd}=a_{\ell,i}^{bd}(y,\DD)$ defined on
$\partial M$ for $0\le i\le\ell$ so that
\begin{eqnarray*}
&&a_n(F,\DD)=\mathcal{I}\{Fa_n\}\text{\ \ and\ \ }
a_{\ell}^{bd}(F,\DD)=\sum_{i=0}^\ell\mathcal{I}^{bd}\{F_ia_{\ell,i}^{bd}\}\,.
\end{eqnarray*}
\end{enumerate}\end{theorem}

We refer to \cite{Gr68,Se69} for a proof of Theorem \ref{thm-1.1} where more general results are obtained in the context
of elliptic operator theory and elliptic boundary conditions. We shall illustrate these formulas in Theorem
\ref{thm-1.3} below. We add a caution that the notation we have chosen differs slightly from what is employed elsewhere.

\subsection{The Bochner Laplacian}\label{sect-1.4}
Before discussing the formulas in Theorem \ref{thm-1.1} in further detail, we must introduce the
formalism of the Bochner Laplacian which will permit us to work in a tensorial and coordinate free fashion.
If $\nabla$
is a connection on
$V$, then we use $\nabla$ and the Levi--Civita connection defined by the metric to covariantly differentiate tensors of all
types. Let `;' denote the components of multiple covariant differentiation -- in particular, $\phi_{;\mu\nu}$ are the components of
$\nabla^2\phi$. If $E$ is an auxiliary endomorphism of $V$, we define the associated {\it modified
Bochner Laplacian} by setting
\begin{equation}\label{eqn-1.e}
\DD(g,\nabla,E)\phi:=-g^{\mu\nu}\phi_{;\nu\mu}-E\phi\,.
\end{equation}

Let $\Gamma_{\mu\nu\sigma}$ and $\Gamma_{\mu\nu}{}^\sigma$ be the Christoffel symbols of the Levi-Civita connection on $M$.
We have, adopting the notation of Equations (\ref{eqn-1.a}) and (\ref{eqn-1.e}), the following (see
\cite{G94}):
\begin{lemma}\label{lem-1.2}
If $\DD$ is an operator of Laplace type, then there exists a unique connection
$\nabla$ on $V$ and a unique endomorphism $E$ on $V$ so that $\DD=\DD(g,\nabla,E)$. The associated connection $1$-form
$\omega$ of $\nabla=\nabla(\DD)$ and the associated endomorphism $E=E(\DD)$ are given by:
\begin{enumerate}
\item $\omega_\mu=\textstyle\frac12(g_{\mu\nu}A_1^\nu+g^{\sigma\varepsilon}\Gamma_{\sigma\varepsilon\mu}
\operatorname{Id})$.
\smallskip\item $E=A_0-g^{\mu\nu}(\partial_{x_\nu}\omega_{\mu}+\omega_{\mu}\omega_{\nu}
    -\omega_{\sigma}\Gamma_{\mu\nu}{}^\sigma)$.
\end{enumerate}\end{lemma}

\subsection{Formulas for the heat trace asymptotics in the smooth setting}\label{sect-1.5}
Let indices $i$, $j$, $k$, $l$ range from $1$ to $m$ and index a local orthonormal
frame $\{e_1,...,e_m\}$ for the tangent bundle of $M$. Let
$R_{ijkl}$ be the components of the Riemann curvature tensor; our sign convention is chosen so that
$R_{1221}=+1$ on the sphere of radius $1$ in $\mathbb{R}^3$. On the collar $\mathcal{C}_\varepsilon$, we normalize the choice of
the local frame by requiring that $e_m=\partial_r$ is the inward unit geodesic normal. We let indices $a$, $b$, $c$, $d$ range
from
$1$ through $m-1$ and index the restricted orthonormal frame $\{e_1,...,e_{m-1}\}$ for the tangent bundle of
the boundary. Let $L_{ab}:=g(\nabla_{e_a}e_b,e_m)$ be the components of the second fundamental form. The scalar invariant $L_{aa}$
is the unnormalized mean curvature (i.e. the geodesic curvature) and will play a central role in our investigation. One has the
following formulas; note that
$F_2=\frac12F_{;mm}$:

\begin{theorem}\label{thm-1.3}
Let $M$ be a compact smooth Riemannian manifold. Let $\DD$ be the Dirichlet realization of an operator of Laplace type. Let $F\in
C^\infty(M)$.
\begin{enumerate}
\smallskip\item $a_0(F,\DD)=(4\pi)^{-m/2}\mathcal{I}\{\Tr(F\Id)\}$.
\smallskip\item $a_1(F,\DD)=\frac16(4\pi)^{-m/2}\mathcal{I}\{\Tr(6FE+FR_{ijji}\Id)\}$.
\smallskip\item $a_0^{bd}(F,\DD)=-\frac14(4\pi)^{-(m-1)/2}\mathcal{I}^{bd}\{\Tr(F_0\Id)\}$.
\smallskip\item $a_1^{bd}(F,\DD)=\frac16(4\pi)^{-m/2}\mathcal{I}^{bd}\{
\Tr(2F_0L_{aa}\Id-3F_1\Id)\}$.
\smallskip\item $a_2^{bd}(F,\DD)=-\frac1{384}(4\pi)^{-(m-1)/2}\mathcal{I}^{bd}\{
\Tr(F_0(96E+[16R_{ijji}-8R_{amma}$
\smallbreak\qquad$+7L_{aa}L_{bb}-10L_{ab}L_{ab}]\Id)-30F_1L_{aa}\Id+48F_2\Id)\}$.
\end{enumerate}
\end{theorem}

Formulas for the invariants $a_n(F,\DD)$ and $a_\ell^{bd}(F,\DD)$ are known for $n,\ell=2,3,4,5$. We refer to \cite{K02} for
further details as the literature is vast and beyond the scope of the present paper to survey.

\subsection{Singular weight functions}\label{sect-1.6}
Fix $\alpha\in\mathbb{C}$. Let $F$ be a smooth function on the interior of $M$ such that $Fr^\alpha\in
C^\infty({\mathcal{C}_{\varepsilon}})$; the parameter $\alpha$ controls the growth (if
$\operatorname{Re}(\alpha)>0$) or decay (if
$\operatorname{Re}(\alpha)<0$) of $F$ near the boundary, assuming that $F r^\alpha$ does not vanish identically on the
boundary. We may expand $F\left|_{\mathcal{C}_{\varepsilon}}\right.$ in a modified Taylor series:
$$
F(y,r)\sim\sum_{i=0}^\infty F_i(y)r^{i-\alpha}\quad\text{where}\quad\left.F_i(y)={\textstyle\frac1{i!}}(\partial_r)^i(r^\alpha F)
\right|_{r=0}\,.
$$
\subsection{Geometry near the boundary}\label{sect-1.7}
The Riemannian measure is in general not product near the boundary, i.e. $dx\ne drdy$, and this plays an important
role in our development.
Let indices
$\sigma,\varrho$ range from $1$ to $m-1$ and index the coordinate frame
$\{\partial_{y_1},...,\partial_{y_{m-1}}\}$ for the tangent bundle of the boundary. One may express the metric on the collar
$\mathcal{C}_\varepsilon$ in the form:
$$ds^2_M=g_{\sigma\varrho}(y,r)dy^\sigma\circ dy^\varrho+dr^2\,.$$
Fix $y_0\in\partial M$ and choose the local coordinates so that $g_{\sigma\varrho}(y_0,0)=\delta_{\sigma\varrho}$. Then we have
that:
\begin{eqnarray}
&&L_{\sigma\varrho}=g(\partial_r,\nabla_{\partial_{x_\sigma}}\partial_{x_\varrho})
    =\Gamma_{\sigma\varrho}{}^m=-\textstyle\frac12\partial_rg_{\sigma\varrho},\nonumber\\
&&g_M(y_0,r)=\sqrt{\det\left\{\Id
+\partial_rg_{\sigma\varrho}(y_0,0)\cdot r+O(r^2)\right\}}
=1-rL_{aa}(y_0)+O(r^2),\label{eqn-1.f}\\
&&\dx=(1-rL_{aa})dr\dy+O(r^2)\,.\nonumber
\end{eqnarray}

\begin{example}\label{exm-1.4}
\rm Let $x_1=\zeta\cos\theta$ and $x_2=\zeta\sin\theta$ be the usual polar coordinates on the unit disk in
$\mathbb{R}^2$. One then has that
$ds^2=d\zeta^2+\zeta^2d\theta^2$ so
$\dx=\zeta d\theta d\zeta$. The geodesic distance to the boundary circle is given by $r=1-\zeta$;
thus $g_{\theta\theta}=(1-r)^2$ and $L_{aa}=1$ so
$\dx=(1-r)dr\dy$.\end{example}

\subsection{Regularization}\label{sect-1.8}
Before discussing the asymptotic expansion of the heat trace in the singular case, we must first discuss regularization; an
analogous regularization was required when discussing the heat content for singular initial temperatures in \cite{BGS08}.
Let $H$ be smooth on the interior of $M$ with $Hr^\alpha\in C^\infty(\mathcal{C}_\varepsilon)$. Then
\begin{eqnarray*}
&&\dx=(1-rL_{aa})dr\dy+O(r^2),\\
&&H\dx=\{H_0r^{-\alpha}+(H_1-H_0L_{aa})r^{1-\alpha}\}dr\dy+O(r^{2-\alpha}dr\dy).
\end{eqnarray*}
For $\operatorname{Re}(\alpha)<3$, define:
\begin{equation}\label{eqn-1.g}
\begin{array}{l}
\Ireg\{H\}
:=\mathcal{I}\{H\dx\}[M-\mathcal{C}_{\varepsilon}]\\
\qquad+\mathcal{I} \left\{H\dx
-\left[H_0r^{-\alpha}+(H_1-H_0L_{aa})r^{1-\alpha}\right]dr\dy
\right\}
[\mathcal{C}_{\varepsilon}]\vphantom{\vrule height 13pt}\\
\qquad+\mathcal{I}^{bd}\{H_0\}
\times\left\{\begin{array}{lll}\frac{\varepsilon^{1-\alpha}}{1-\alpha}&\text{if}&\alpha\ne1,\\
\ln(\varepsilon)&\text{if}&\alpha=1.
\end{array}\right.\\
\qquad+\mathcal{I}^{bd}\{H_1-H_0L_{aa}\}\times\left\{\begin{array}{lll}
\frac{\varepsilon^{2-\alpha}}{2-\alpha}&\text{if}&\alpha\ne2,\\
\ln(\varepsilon)&\text{if}&\alpha=2.
\end{array}\right.
\end{array}\end{equation}
This is independent of the parameter $\varepsilon$ and agrees with
$\mathcal{I}\{H\}$ if $\operatorname{Re}(\alpha)<1$. Because the
integrand over $\mathcal{C}_\varepsilon$ is
$O(r^{2-\operatorname{Re}(\alpha)})$ and
$\operatorname{Re}(2-\alpha)>-1$, $\Ireg$ is well defined.

The regularization $\Ireg$ is a
meromorphic function of $\alpha$ with simple poles at $\alpha = 1,2$. At these exceptional values,
$\Ireg$ is defined as the constant term in the appropriate Laurent expansion, thus
dropping the pole. We shall apply this regularization to functions of the form $H(x)=F(x)a_n(x,\DD)$.

\subsection{Heat trace asymptotics in the singular setting}\label{sect-1.9}
The Dirichlet heat kernel satisfies $p_\DD(x,(\tilde y,\tilde r),t)|_{\tilde r=0}=0$. Since $p_{\mathcal{D}}$ is smooth for $t>0$
and $\mathcal{C}_{\varepsilon}$ is compact, we may use the Taylor
series expansion of $p_{\mathcal{D}}$ to derive the estimate:
$$|p_\DD(x,(\tilde y,\tilde r);t)|\le C(t)\tilde r\quad\text{on}\quad\mathcal{C}_\varepsilon\,.$$
A similar estimate holds for $|p_\DD((y,r),\tilde x;t)|$. We set
$\tilde x=(y,r)$ to derive the estimate on the diagonal:
\begin{equation}\label{eqn-1.h}
|p_\DD((y,r),(y,r);t)|\le C(t)r^2\quad\text{on}\quad\mathcal{C}_\varepsilon\,.
\end{equation}
Thus if $\operatorname{Re}(\alpha)<3$, then
Equation (\ref{eqn-1.c}) shows that $\Tr_{L^2}(Fe^{-t\DD})$
is convergent.

We shall begin our investigation in Section \ref{sect-2} by
establishing the following result by a direct computation as this
motivates our entire investigation; note that only limited
smoothness is required of the boundary in this result. It has the
further advantage of confirming by completely different means some
of the constants that will be computed again in Sections 4 and 5.
Let $C$ here and elsewhere denote Euler's constant.

\begin{theorem}\label{thm-1.5}  Let $M\subset\mathbb{R}^2$ be an open, bounded, and connected planar set with
$C^2$ boundary. Let
$0<\varepsilon_0<\varepsilon$. Set
$F(x):=F_0(y)r^{-\alpha}\chi(r)$ where
$\chi\in C^\infty(\mathcal{C}_\varepsilon)$ satisfies:
$$\chi(r)=\left\{\begin{array}{lll}
1&\text{if}&0\le r\le\varepsilon_0,\\0&\text{if}&\varepsilon\le r\,.\end{array}\right.$$
\begin{enumerate}
\item If $0<\alpha<1$ and $t\downarrow0$, then:
\begin{eqnarray*}
&&\textstyle\Tr_{L^2}(Fe^{-t\Delta_M})=\frac1{4\pi t}\{
\mathcal{I}\{F\}-\frac12\Gamma\left(\frac{1-\alpha}2\right)t^{(1-\alpha)/2}{\mathcal{I}}^{bd}\{F_0\}\\
&&\quad\textstyle+\frac{4-\alpha}{4(3-\alpha)}\Gamma\left(\frac{2-\alpha}2\right)t^{1-\alpha/2}
{\mathcal{I}}^{bd}\{F_0L_{aa}\}\}+O(1)\,.
\end{eqnarray*}
\item If $\alpha=1$
and $t\downarrow0$, then:
$$
\textstyle
\Tr_{L^2}\left(Fe^{-t\Delta_M}\right)=\frac1{4\pi
t}\left\{\Ireg\{F\}-{\textstyle\frac12}\ln(t)\cdot{\mathcal{I}}^{bd}\{F_0\}+{\textstyle\frac C2}{\mathcal{I}}^{bd}\{F_0\}
\right\}+O(t^{-\frac12})\,.
$$
\end{enumerate}\end{theorem}

This result extends to a very general setting. The following Theorem generalizes
Theorem \ref{thm-1.1} to the singular setting where, in contrast to Theorem \ref{thm-1.5}, we assume the boundary is $C^\infty$.
We also refer to
\cite{BGS08} for further details where an analogous result was proved for the heat content asymptotics.
In Section \ref{sect-5}, we will use the pseudo-differential calculus to show that:

\begin{theorem}\label{thm-1.6}
Let $M$ be a compact smooth Riemannian manifold. Let $\DD$ be the Dirichlet realization of an operator of Laplace type. Let
$a_n=a_n(x,\DD)$ be the interior local heat trace asymptotics of Theorem 1.1.  Fix
$\alpha\in\mathbb{C}$ with $\operatorname{Re}(\alpha)<3$. Let $F$ be a smooth function on the interior of $M$ such that
$Fr^\alpha\in C^\infty({\mathcal{C}_{\varepsilon}})$.
\begin{enumerate}
\item If $\alpha\ne1,2$, there is a complete asymptotic
expansion as
$t\downarrow0$ of the form:\par\noindent
\begin{eqnarray*}&&\Tr_{L^2}(Fe^{-t\DD})\sim t^{-m/2}\sum_{n=0}^\infty t^n\mathcal{I}_{\operatorname{Reg}}
\{Fa_{n}\}
+t^{-(m-1)/2}\sum_{\ell=0}^\infty t^{(\ell-\alpha)/2}a_{\ell,\alpha}^{bd}(F,\DD)\,.
\end{eqnarray*}
\item If $\alpha=1,2$, there is a complete asymptotic
expansion as
$t\downarrow0$ of the form:
\begin{eqnarray*}&&\Tr_{L^2}(Fe^{-t\DD})\sim t^{-m/2}\sum_{n=0}^\infty t^n\mathcal{I}_{\operatorname{Reg}}
\{Fa_{n}\}
+t^{-(m-1)/2}\sum_{\ell=0}^\infty t^{(\ell-\alpha)/2}a_{\ell,\alpha}^{bd}(F,\DD)\\
&&\qquad\qquad\qquad+t^{-m/2}\ln(t)\sum_{k=0}^\infty t^{k/2}\tilde a_{k,\alpha}^{bd}(F,\DD)\,.
\end{eqnarray*}
\item There exist local invariants $a_{\ell,\alpha,i}^{bd}=a_{\ell,\alpha,i}^{bd}(y,\DD)$ on $\partial M$,
which are holomorphic in $\alpha$ for
$\alpha\ne1,2$, so that
$$a_{\ell,\alpha}^{bd}(F,\DD)=\displaystyle\sum_{i=0}^\ell
\mathcal{I}^{bd}\{F_ia_{\ell,\alpha,i}^{bd}\}\,.$$ The invariants
$a_{\ell,z,i}^{bd}$ have simple poles at $z=1,2$ and
$$a_{\ell,\alpha,i}^{bd}
=\left.\left\{a_{\ell,z,i}^{bd}-\frac1{z-\alpha}\operatorname{Res}_{z=\alpha}a_{\ell,z,i}^{bd}\right\}\right|_{z=\alpha}
\quad\text{if}\quad\alpha=1,2 \,.$$
\item The $\ln(t)$ coefficients in Assertion (2) are given by:
\begin{eqnarray*}
&&\tilde a_{k,\alpha}^{bd}(F,\DD)=\left\{\begin{array}{rl}
-{\textstyle\frac12}\mathcal{I}^{bd}\{(Fa_n)_0\}&\text{if }k=2n\text{ and }\alpha=1,\\
-{\textstyle\frac12}\mathcal{I}^{bd}\{(Fa_n)_1-(Fa_n)_0L_{aa}\}&\text{if }k=2n\text{ and }\alpha=2,\\
 0&\text{if }k=2n+1.\end{array}\right.
\end{eqnarray*}
\end{enumerate}\end{theorem}

Throughout this paper, let
$$\textstyle\kappa_\alpha:=\frac12\Gamma\left(\frac{1-\alpha}2\right)\,.$$
The boundary invariants for $\alpha\ne1,2$ and
for $\ell=0,1,2$ are given by:

\begin{theorem}\label{thm-1.7}
If $\alpha\ne1,2$, then one has:
\begin{enumerate}
\item $a_{0,\alpha}^{bd}(F,\DD)=\kappa_\alpha(4\pi)^{-m/2}
\mathcal{I}^{bd}\{\Tr(-F_0\Id)\}$.
\smallbreak\item $a_{1,\alpha}^{bd}(F,\DD)=\kappa_{\alpha-1}(4\pi)^{-m/2}\mathcal{I}^{bd}
\{\Tr(-F_1\Id+\frac{\alpha-4}{2(\alpha-3)}F_0L_{aa}\Id)\}$.
\smallbreak\item $a_{2,\alpha}^{bd}(F,\DD)=\kappa_{\alpha-2}(4\pi)^{-m/2}\mathcal{I}^{bd}\{\Tr(
-F_2\Id+\frac{\alpha-5}{2(\alpha-4)}F_1L_{aa}\Id
$\smallbreak\quad$
+\frac{1}{6}F_0R_{amma}\Id
-\frac{\alpha-7}{8(\alpha-6)}F_0L_{aa}L_{bb}\Id
+\frac{\alpha-5}{4(\alpha-6)}F_0L_{ab}L_{ab}\Id$\smallbreak\quad$
-\frac1{3(1-\alpha)}F_0R_{ijji}\Id-\frac2{1-\alpha}F_0E)\}$.
\end{enumerate}
\end{theorem}
We remark that we recover Theorem \ref{thm-1.3} by setting
$\alpha=0$ in Theorem \ref{thm-1.7}. We omit details as the
calculation is entirely elementary. The boundary invariants for
$\alpha=1,2$, and for $\ell=0,1,2$ are given by:
\begin{theorem}\label{thm-1.8}
\ \begin{enumerate}\item When $\alpha=1$,
\begin{enumerate}
\smallbreak\item $a_{0,1}^{bd}(F,\DD)=(4\pi)^{-m/2}\mathcal{I}^{bd}\{\Tr(\frac C2 F_0\Id)\}$.
\smallbreak\item $a_{1,1}^{bd}(F,\DD)=(4\pi)^{-m/2}\frac{\sqrt\pi}2\mathcal{I}^{bd}\{
   \Tr(-F_1\Id+\frac34F_0L_{aa}\Id)\}$.
\smallbreak\item $a_{2,1}^{bd}(F,\DD)=(4\pi)^{-m/2}\mathcal{I}^{bd}\left\{
\Tr\left(\textstyle-\frac12F_2\Id+\frac13F_1L_{aa}\Id+\frac1{12}R_{amma}\Id\right.\right.\\
\qquad\qquad\qquad\qquad\textstyle\left.\left.
-\frac3{40}L_{aa}L_{bb}\Id+\frac1{10}L_{ab}L_{ab}\Id+{\textstyle\frac C{12}} R_{ijji}\Id+{\textstyle\frac C2}
E\right)\right\}\,.$
\end{enumerate}
\smallbreak\item When $\alpha=2$,
\begin{enumerate}
\smallbreak\item $a_{0,2}^{bd}(F,\DD)=(4\pi)^{-m/2}\sqrt\pi\mathcal{I}^{bd}\{\Tr(F_0\Id)\}$.
\smallbreak\item $a_{1,2}^{bd}(F,\DD)=(4\pi)^{-m/2}\mathcal{I}^{bd}\{\Tr({\textstyle\frac C2} F_1\Id-[{\textstyle\frac
C2}+\frac12]F_0L_{aa}\Id)\}$.
\smallbreak\item $a_{2,2}^{bd}(F,\DD)=(4\pi)^{-m/2}\sqrt\pi\mathcal{I}^{bd}\{\Tr(-\frac12F_2\Id+\frac38F_1L_{aa}\Id$
\smallbreak\quad$+F_0(\frac1{12}R_{amma}-\frac5{64}L_{aa}L_{bb}+\frac3{32}L_{ab}L_{ab}+\frac16R_{ijji})\Id+F_0E)\}$.
\end{enumerate}\end{enumerate}
\end{theorem}

Here is a brief guide to the remainder of this paper.
In Section \ref{sect-2}, we will make a special case calculation to establish Theorem \ref{thm-1.5}. A probabilistic estimate
of R.\ Lang  \cite{RL11} and of  H.\ R.\ Lerche and D. Siegmund
\cite{HLDS10} plays a central role.
 In Section \ref{sect-3}, we shall use dimensional analysis (scaling
arguments) and various functorial properties to study the heat trace invariants. We will derive Theorem \ref{thm-1.6} (4) from
the asymptotic series in Theorem
\ref{thm-1.6} (3); another derivation will be given subsequently in Section \ref{sect-5}. We shall examine the general form of
the invariants and establish the following result.

\begin{lemma}\label{lem-1.9}
\ \begin{enumerate}
\item There exist
universal constants $\{\bar\kappa_\alpha,\kappa_\alpha^1,\kappa_\alpha^3,\kappa_\alpha^4,\kappa_\alpha^5\}$ so that:
\begin{enumerate}
\item $a_{0,\alpha}^{bd}(F,\DD)=(4\pi)^{-m/2}\mathcal{I}^{bd}\{\Tr(-\bar\kappa_\alpha F_0\Id)\}$.
\smallbreak\item $a_{1,\alpha}^{bd}(F,\DD)=(4\pi)^{-m/2}\mathcal{I}^{bd}
\{\Tr(-\bar\kappa_{\alpha-1}F_1\Id+\kappa_\alpha^1F_0L_{aa}\Id)\}$.
\smallbreak\item
$a_{2,\alpha}^{bd}(F,\DD)=(4\pi)^{-m/2}\mathcal{I}^{bd}\{\Tr(-\bar\kappa_{\alpha-2}F_2\Id+\kappa_{\alpha-1}^1F_1L_{aa}\Id$
\smallbreak$
+F_0[\kappa_\alpha^3R_{amma}+\kappa_\alpha^4L_{aa}L_{bb}
+\kappa_\alpha^5L_{ab}L_{ab}]\Id-\bar\kappa_\alpha F_0[E+\frac16R_{ijji}\Id])\}$.
\end{enumerate}
\item If $\alpha\ne1,2$, then $\bar\kappa_\alpha=\kappa_\alpha$ and
$\kappa_{\alpha}^1=\frac12\Gamma\left(\frac{2-\alpha}2\right)\frac{\alpha-4}{2(\alpha-3)}$.
\item $\bar\kappa_1=\frac C2$.
\end{enumerate}
\end{lemma}


In Section \ref{sect-4}, we evaluate the remaining
universal coefficients of Lemma \ref{lem-1.9} using the calculus of
pseudo-differential operators and complete the proof of Theorem
\ref{thm-1.7} by showing:
\begin{lemma}\label{lem-1.10}
Adopt the notation of Lemma \ref{lem-1.9}. If $\alpha\ne1,2$, then:
\begin{eqnarray*}
&&\textstyle\kappa_\alpha^3=-\frac{\alpha-1}{24}\Gamma\left(\frac{1-\alpha}2\right),\quad
  \kappa_\alpha^4=\frac{7-8\alpha+\alpha^2}{32(\alpha-6)}\Gamma\left(\frac{1-\alpha}2\right),\quad
\textstyle\kappa_\alpha^5=\frac{6\alpha-5-\alpha^2}{16(\alpha-6)}\Gamma\left(\frac{1-\alpha}2\right)\,.
\end{eqnarray*}
\end{lemma}
\noindent We conclude the paper in Section \ref{sect-5} by using the pseudo-differential calculus to establish Theorem
\ref{thm-1.6}. We have postponed the proof of Theorem \ref{thm-1.6} until this point as much of the needed notation will be
established in Section \ref{sect-4}. We will also complete the proof of Theorem \ref{thm-1.8}.

We have chosen to use special case calculations, the functorial method, and the pseudo-differential calculus as our
purpose in this paper is at least in part expository and we wish to illustrate the interplay amongst these methods.
In a subsequent paper, we shall perform a similar analysis for other elliptic boundary conditions (Robin, transfer, transmittal,
etc.); it will be necessary to restrict to
$\operatorname{Re}(\alpha)<1$ to ensure convergence and regularization will not be required in that analysis.

\section{Computations in $\mathbb{R}^2$}\label{sect-2}
This section is devoted to the proof of Theorem \ref{thm-1.5}, and we shall adopt the notation of that Theorem throughout. As we
shall be dealing with different weights, we drop the notation
$\mathcal{I}$ and return to ordinary integrals in this section to perform a special case calculation in flat space.
One has the following estimate of R.\ Lang  \cite{RL11} and of  H.\ R.\ Lerche and D. Siegmund
\cite{HLDS10} that adjusts the formula of Equation (\ref{eqn-1.d}) for the heat kernel on a half space to take into account the
curvature of the boundary of $M$ (for related results see also
\cite{BG07,L68}):

\begin{lemma}\label{lem-2.1}
Adopt the notation of Theorem \ref{thm-1.5}. Let $x=(y,r)\in\mathcal{C}_\varepsilon$. As $t\downarrow0$,
$$
p_M(x,x;t)=\frac1{4\pi t}
\left\{1-e^{-r^2/t}-L_{aa}(y)r^2t^{-1/2}\int^\infty_{rt^{-1/2}}
e^{-\eta^2} d \eta\right\}+O(1)\,.
$$
\end{lemma}

\medbreak\noindent{\it Proof of Theorem \ref{thm-1.5} (1).} Parametrize the
boundary of
$M$ by arclength. There is no higher order correction in $\mathbb{R}^2$ and the $O(r^2)$ term in Equation (\ref{eqn-1.f})
vanishes. Thus on the collar
$\mathcal{C}_\varepsilon$, we have
$$dx=(1-L_{aa}(y)r)dr\dy\,.$$
Following Equation (\ref{eqn-1.d}), we let $p_H(x,x;t)=(4\pi t)^{-1}(1-e^{-r^2/t})$ be the Dirichlet heat kernel in
the halfspace $r\ge0$ on the diagonal. We take $0<\alpha<1$ and express:
\begin{eqnarray*}
&&\Tr_{L^2}(Fe^{-t\Delta_M})\\
&&\qquad=\int_{\partial M}\int_0^\varepsilon F_0(y)r^{-\alpha}\chi(r)(1-L_{aa}(y)r)p_M((y,r),(y,r);t)
dr\dy\\
&&\qquad=\int_{\partial M}\int_0^\infty F_0(y)r^{-\alpha}\chi(r)(1-L_{aa}(y)r)p_M((y,r),(y,r);t)
dr\dy\\
&&\qquad=D_1+D_2+D_3+D_4+D_5
\end{eqnarray*}
where, motivated by Lemma \ref{lem-2.1}, we have:
\medbreak\qquad
$\displaystyle D_1:=\frac1{4\pi t}\int_MF(x)dx$,
\medbreak\qquad
$\displaystyle D_2:=-\frac1{4\pi t}\int_{\partial M}\int_0^\infty F_0(y) r^{-\alpha}e^{-r^2/t}dr\dy$,
\medbreak\qquad
$\displaystyle D_3:=\frac1{4\pi t}\int_{\partial M}\int_0^\infty F_0(y)r^{-\alpha}(1-\chi(r))(1-L_{aa}(y)r)e^{-r^2/t}dr\dy$,
\medbreak\qquad
$\displaystyle D_4:=\frac1{4\pi t}\int_{\partial M}\int_0^\infty F_0(y)L_{aa}(y)r^{1-\alpha}e^{-r^2/t}dr\dy$,
\medbreak\qquad
$\displaystyle D_5:=\int_{\partial M}\int_0^\infty F_0(y)r^{-\alpha}\chi(r)(1-L_{aa}(y)r)
\{(p_M-p_H)((y,r),(y,r);t)\}drdy$.
\medbreak\noindent
A straightforward computation yields:
\medbreak\qquad
$\displaystyle D_2=-\frac1{4\pi}\cdot\frac12\Gamma\left(\frac{1-\alpha}2\right)t^{-(1+\alpha)/2}\int_{\partial M}F_0(y)\dy$,
\medbreak\qquad
$\displaystyle D_3=O\left(e^{-\varepsilon_0^2/(2t)}\right)$,
\medbreak\qquad
$\displaystyle D_4=\frac1{4\pi}\cdot\frac12\Gamma\left(\frac{2-\alpha}2\right)t^{-\alpha/2}\int_{\partial M}F_0(y)L_{aa}(y)\dy$.
\medbreak\noindent We use Lemma \ref{lem-2.1} to compute $D_5$. Since $F_0(y)r^{-\alpha}\chi(r)$ is integrable
on $\mathcal{C}_\varepsilon$, we have the $O(1)$ in Lemma \ref{lem-2.1} remains $O(1)$ as $t\downarrow 0$. Hence
for $0<\alpha<1$,
\medbreak\qquad
$\displaystyle D_5=-\frac1{4\pi t}\int_{\partial M}\int_0^\infty\int_{rt^{-1/2}}^\infty F_0(y)L_{aa}(y)r^{2-\alpha}t^{-1/2}
e^{-\eta^2}d\eta dr\dy$
\medbreak\qquad\quad
$\displaystyle+\frac1{4\pi t}\int_{\partial
M}\int_0^\infty\int_{rt^{-1/2}}^\infty F_0(y)L_{aa}(y)L_{bb}(y)r^{3-\alpha}t^{-1/2}e^{-\eta^2}d\eta dr\dy+O(1)$
\medbreak\qquad\quad
$\displaystyle=\frac1{4\pi}\cdot\frac{2-\alpha}{4(\alpha-3)}\Gamma\left(\frac{2-\alpha}2\right)
t^{-\alpha/2}\int_{\partial M}F_0(y)L_{aa}(y)dy+O(t^{(1-\alpha)/2})+O(1)$
\medbreak\qquad\quad
$\displaystyle=\frac1{4\pi}\cdot\frac{2-\alpha}{4(\alpha-3)}\Gamma\left(\frac{2-\alpha}2\right)
t^{-\alpha/2}\int_{\partial M}F_0(y)L_{aa}(y)dy+O(1)$.
\medbreak\noindent We collect terms to complete the proof of Theorem \ref{thm-1.5} (1).\hfill\qedbox

\begin{remark}\label{rmk-2.2}
\rm We have chosen to study the region $0<\alpha<1$. The reason for this is that $F$ has to be
integrable in order to control the $O(1)$ remainder in Lemma \ref{lem-2.1}. If one wishes to obtain just the leading asymptotic
behaviour of
$Tr_{L^2}(Fe^{-t\Delta_M})$, then probabilistic estimates for
\begin{equation}\label{eqn-2.a}
R_M(x;t):=(p_M-p_H)((y,r),(y,r);t)
\end{equation}
along the lines of \cite{B-87}, and analogous to \cite{B-07}, could be used to show that for $1\le\alpha<3$,
\begin{equation}\label{eqn-2.b}
\int_{\partial M}\int_0^\varepsilon\chi(r)r^{-\alpha}R(x;t)drdy=O(t^{-(m-2+\alpha)/2}),\quad t\downarrow0\,.
\end{equation}\end{remark}

\medbreak\noindent{\it Proof of Theorem \ref{thm-1.5} (2).} Recall from Equation (\ref{eqn-1.d}) the formula for the heat kernel
$p_H$ on the diagonal for the half space. We decompose
$$\Tr_{L^2}(Fe^{-t\Delta_M})=E_1+E_2+E_3$$
where we have:
\medbreak\qquad
$\displaystyle E_3:=\int_{\partial M}\int_0^\infty F_0(y)r^{-1}\chi(r)(1-L_{aa}(y)r)(p_M-p_H)((y,r),(y,r);t))drdy$,
\medbreak\qquad
$\displaystyle E_2:=-\frac1{4\pi t}\int_{\partial M}\int_0^\infty F_0(y)L_{aa}(y)\chi(r)(1-e^{-r^2/t})drdy$,
\medbreak\qquad
$\displaystyle E_1:=\frac1{4\pi t}\int_{\partial M}\int_0^\infty F_0(y)r^{-1}\chi(r)(1-e^{-r^2/t})drdy$.
\medbreak\noindent
We apply Equations (\ref{eqn-2.a}) and (\ref{eqn-2.b}) with $m=2$ and with $\alpha=1$ to see that $E_3=O(t^{-1/2})$.
Furthermore, $E_2=O(t^{-1/2})$. The leading term is provided by $E_1$. We
integrate by parts to see that:
\medbreak\qquad
$\displaystyle E_1=-\frac1{4\pi t}\int_{\partial M}F_0(y)\int_0^\infty
\ln(r)\left\{\chi^\prime(r)(1-e^{-r^2/t})+\chi(r)\frac{2r}te^{-r^2/t}
\right\}drdy$
\medbreak
$\displaystyle\qquad=-\frac1{4\pi t}\int_{\partial M}F_0(y)\left\{\int_0^\infty\ln(r)\chi^\prime(r)dr+\int_0^\infty
\ln(r)\frac{2r}te^{-r^2/t}dr\right\}dy$\medbreak\qquad\qquad$+O(e^{-\varepsilon_0^2/(2t)})$
\medbreak
$\displaystyle\qquad=\frac1{4\pi t}\cdot\frac12\int_{\partial
M}F_0(y)\left\{\ln\left({\textstyle\frac{\varepsilon_0^2}t}\right)+C+2\int_{\varepsilon_0}^\varepsilon
\frac{\chi(r)}rdr\right\} +O(e^{-\varepsilon_0^2/(2t)})$.
\medbreak\noindent This completes the proof of Theorem \ref{thm-1.5}.\hfill\qedbox

\section{The functorial method}\label{sect-3}
We adopt the notation of Theorem \ref{thm-1.6} throughout this section.
We begin our study with the following:
\begin{lemma}\label{lem-3.1}
There exist constants $\varepsilon_{\ell,\alpha}^\nu$ so that
\begin{enumerate}
\item $a_{0,\alpha}^{bd}(F,\DD)=(4\pi)^{-m/2}\mathcal{I}^{bd}\{\Tr(\varepsilon_{0,\alpha}^0F_0\Id)\}$.
\smallbreak\item
$a_{1,\alpha}^{bd}(F,\DD)=(4\pi)^{-m/2}\mathcal{I}^{bd}\{\Tr(\varepsilon_{1,\alpha}^0F_1\Id
   +\varepsilon_{1,\alpha}^1F_0L_{aa}\Id)\}$.
\smallbreak\item $a_{2,\alpha}^{bd}(F,\DD)=(4\pi)^{-m/2}\mathcal{I}^{bd}\{\Tr(
\varepsilon_{2,\alpha}^0F_2\Id+\varepsilon_{2,\alpha}^1F_1L_{aa}\Id$\smallbreak$
+F_0[\varepsilon_{2,\alpha}^2R_{ijji}+\varepsilon_{2,\alpha}^3R_{amma}+\varepsilon_{2,\alpha}^4L_{aa}L_{bb}
+\varepsilon_{2,\alpha}^5L_{ab}L_{ab}]\Id
+\varepsilon_{2,\alpha}^6F_0E)\}$.
\end{enumerate}\end{lemma}

\begin{proof} We apply dimensional analysis -- we shall suppose that $\alpha\notin\mathbb{Z}$ for the moment.  Let $c>0$ define a
rescaling
$g_c:=c^2g$. We then have
$$
\begin{array}{lll}
\dx_c=c^m\dx,&\dy_c=c^{m-1}\dy,&\DD_c=c^{-2}\DD,\\
r_c:=cr,&\partial_{r_c}=c^{-1}\partial_r,&F_{i,c}=c^{\alpha-i}F_i,\\
\Iregc=c^m\Ireg,&\mathcal{I}_c^{bd}=c^{m-1}\mathcal{I}^{bd}\,.
\end{array}$$
Let $a_{n,c}:=a_n(x,\DD_c)$, $a_{\ell,\alpha,c}:=a_{\ell,\alpha}(y,\DD_c)$, and $a_{\ell,\alpha,i,c}:=a_{\ell,\alpha,i}(y,\DD_c)$
denote the local heat trace invariants defined by
$\DD_c$ on
$M$ and on
$\partial M$, respectively. It is immediate that
\begin{equation}\label{eqn-3.a}
\Tr_{L^2}(Fe^{-(tc^{-2})\DD})=\Tr_{L^2}(Fe^{-t\DD_c})\,.
\end{equation}
We expand both sides of Equation (\ref{eqn-3.a}) in an asymptotic expansion:
\medbreak\qquad
$\displaystyle\phantom{\sim} t^{-m/2}c^m\sum_{n=0}^\infty(c^{-2}t)^n\Ireg\{Fa_n\}$
\medbreak\qquad$\displaystyle
\qquad+(c^{-2}t)^{-(m-1)/2}\sum_{\ell=0}^\infty c^{\alpha-\ell}t^{(\ell-\alpha)/2}\sum_{i=0}^\ell\mathcal{I}^{bd}
\{F_ia_{\ell,\alpha,i}\}$
\smallbreak\qquad$\displaystyle
\sim t^{-m/2}\sum_{n=0}^\infty t^nc^m\Ireg\{Fa_{n,c}\}$
\smallbreak\qquad$\displaystyle
\qquad+ t^{-(m-1)/2}\sum_{\ell=0}^\infty t^{(\ell-\alpha)/2}\sum_{i=0}^\ell c^{m-1}c^{\alpha-i}
\mathcal{I}^{bd}\{F_{i}a_{\ell,\alpha,i,c}\}$.\medbreak\noindent
Since $\alpha\notin\mathbb{Z}$, the interior and the boundary terms decouple. We equate terms in the asymptotic expansions
to see that
$$a_{n,c}=c^{-2n}a_n\quad\text{and}\quad a_{\ell,\alpha,i,c}=c^{i-\ell}a_{\ell,\alpha,i}\,.$$
Examining relations of this kind is straightforward -- they mean that the local formula $a_n(x,\DD)$ is homogeneous of weighted
degree
$2n$ in the jets of the derivatives of the symbol of $\DD$ and that the local formula $a_{\ell,\alpha,i}(y,\DD)$ is homogeneous
of weighted degree
$\ell-i$ in the jets of the derivatives of the symbol of $\DD$. One may use Weyl's theory of invariants to express a spanning
set for the invariants which arise in this context and complete the proof of Lemma \ref{lem-3.1} for $\alpha\notin\mathbb{Z}$.
We use analytic continuation to establish Lemma \ref{lem-3.1} when
$\alpha=0,-1,-2,...$ as well. We refer to \cite{G94} for further details concerning this sort of dimensional analysis.

If $\alpha=1$, then the argument is rather different. Let $\varepsilon_c$ be the width of the collar $\mathcal{C}$ with respect
to the rescaled metric. The regularizing term in Equation (\ref{eqn-1.g}) does not simply rescale. Rather we have:
\begin{eqnarray*}
&&\ln(\varepsilon_c)\mathcal{I}_c^{bd}\{(Fa_{n,c})_{0,c}\}
  =\ln(c\varepsilon)c^\alpha c^{-2n}c^{m-1}\mathcal{I}^{bd}\{(Fa_n)_0\}\\
&=&c^{m-2n}\{\ln(c)+\ln(\varepsilon)\}\mathcal{I}^{bd}\{(Fa_n)_0\}\,.
\end{eqnarray*}
This yields the modified relation:
$$
\Iregc\{Fa_{n,c}\}=c^{m-2n}\Ireg\{Fa_n\}+\ln(c)c^{m-2n}\mathcal{I}^{bd}\{(Fa_n)_0\}\,.
$$
A similar argument for $\alpha=2$ shows that:
$$
\Iregc\{Fa_{n,c}\}=c^{m-2n}\Ireg\{Fa_n\}+\ln(c)c^{m-2n}\mathcal{I}^{bd}\{(Fa_n)_1-(Fa_n)_0L_{aa}\}\,.
$$
When we use Equation (\ref{eqn-3.a}) to equate coefficients in the asymptotic series, we have
$$a_{\ell,\alpha}^{bd}(F,\DD_c)=c^\ell a_{\ell,\alpha}^{bd}(F,\DD),$$ which
completes the proof of Lemma
\ref{lem-3.1} in these exceptional cases. We compare the terms involving $\ln(c)$ to obtain
additional relations. The left hand side in the following equation arises from $\ln(c^{-2}t)$ and the right hand side arises
from $\mathcal{I}_{reg,c}$ when we apply Equation (\ref{eqn-3.a}); if $k=2n$, then
\begin{eqnarray*}
&&-2\ln(c)\tilde a_{k,1}^{bd}(F,\DD)=\ln(c)\mathcal{I}^{bd}\{(Fa_n)_0\}\\
&&-2\ln(c)\tilde a_{k,2}^{bd}(F,\DD)=\ln(c)\mathcal{I}^{bd}\{(Fa_n)_1-(Fa_n)_0L_{aa}\}\,.
\end{eqnarray*}
There are no corresponding terms if $k$ is odd and thus $\tilde
a_{k,1}=0$  and $\tilde a_{k,2}=0$ if $k$ is odd. This establishes
Theorem \ref{thm-1.6} (4). An alternate proof is given in Section
5.\end{proof}

 We now use the functorial method to establish the following
result:
\begin{lemma}\label{lem-3.2}
\ \begin{enumerate}
\item $\varepsilon_{2,\alpha}^2=\frac16\varepsilon_{0,\alpha}^0$ and $\varepsilon_{2,\alpha}^6=\varepsilon_{0,\alpha}^0$.
\item The constants $\varepsilon_{\ell,\alpha}^\mu$ of Lemma \ref{lem-3.1} are dimension free.
\item $\varepsilon_{1,\alpha}^0=\varepsilon_{0,\alpha-1}^0$, $\varepsilon_{2,\alpha}^0=\varepsilon_{0,\alpha-2}^0$, and
$\varepsilon_{2,\alpha}^1=\varepsilon_{1,\alpha-1}^1$.
\end{enumerate}
\end{lemma}

\begin{proof} Suppose that $M=M_1\times M_2$, that $g_M=g_{M_1}+g_{M_2}$, that $\DD_M=\DD_{M_1}+\DD_{M_2}$, and that $F_M=F_1F_2$
where
$F_i$ are defined on $M_i$. We suppose that
$M_1$ is a closed manifold and thus $\partial M=M_1\times\partial M_2$.
We then have:
\begin{eqnarray*}
&&\displaystyle e^{-t\DD_M}=e^{-t\DD_{M_1}}e^{-t\DD_{M_2}},\\
&&\displaystyle\Tr_{L^2}(F_Me^{-t\DD_M})=\Tr_{L^2}(F_1e^{-t\DD_{M_1}})\cdot\Tr_{L^2}(F_2e^{-t\DD_{M_2}})\,.
\end{eqnarray*}
Equating asymptotic series yields
$$\displaystyle a_{\ell,\alpha}^{bd}(F_M,\DD_M)=\sum_{2k+j=\ell}a_k(F_1,\DD_{M_1})a_{j,\alpha}^{bd}(F_2,\DD_{M_2})$$
and hence a corresponding decomposition of the local formulas:
\begin{equation}\label{eqn-3.b}
\displaystyle a_{\ell,\alpha,i}^{bd}(y,\DD_M)=\sum_{2k+j=\ell}a_k(x_1,\DD_{M_1})a_{j,\alpha,i}^{bd}(y_2,\DD_{M_2})
\quad\text{for}\quad y=(x_1,y_2)\,.
\end{equation}
Assertion (1) now follows from Theorem \ref{thm-1.3} (2), from Lemma \ref{lem-3.1}, and from Equation (\ref{eqn-3.b});
the multiplicative constants $(4\pi)^{-m/2}$ play no role.
If we take $M_1=S^1$ and $\DD_{M_1}=-\partial_\theta^2$, then the structures are flat. Thus
$a_0(x,\DD_{M_1})=1/\sqrt{4\pi}$ and $a_k(x,\DD_{M_1})=0$ for $k\ge1$. Thus Equation (\ref{eqn-3.b}) yields in this special case
the following identity from which Assertion (2) follows after taking into account the multiplicative constants $(4\pi)^{-m/2}$:
$$a_{\ell,\alpha,i}^{bd}(y,\DD_M)=\textstyle\frac1{\sqrt{4\pi}}a_{\ell,\alpha,i}^{bd}(y_2,\DD_{M_2})\,.$$

We prove Assertion (3) by index shifting. Let $\chi(r)$ be a smooth function so that $\chi\equiv0$ near $r=\varepsilon$. Let
$F(y,r)=F_0(y)\chi(r)r^{-\alpha_1}$ for $\operatorname{Re}(\alpha_1)<3$. We apply Theorem \ref{thm-1.6} with $\alpha=\alpha_1$
and with
$\alpha=\alpha_1-1$ to see that:
$$a_{\ell,\alpha,j}^{bd}(y,\DD)=a_{\ell-1,\alpha-1,j-1}^{bd}(y,\DD)\quad\text{for}\quad j\ge1\,.$$
Assertion (3) now follows.
\end{proof}

\medbreak\noindent{\it Proof of Lemma \ref{lem-1.9}.} Assertion (1) of Lemma \ref{lem-1.9} follows from Lemma \ref{lem-3.1}
and Lemma
\ref{lem-3.2} by a suitable relabeling of the coefficients. Assertions (2) and (3) follow from Theorem \ref{thm-1.5}.
\hfill\qedbox

\medbreak \noindent{\it Proof of Theorem \ref{thm-1.8}.} We now
derive Theorem \ref{thm-1.8} from Theorem \ref{thm-1.7} using
Assertion (3) of Theorem 1.6. Certain of the coefficients are
regular so the computation is elementary; we do not need to drop the
pole. We simply set $\alpha=1$ to compute $a_{1,1}^{bd}$ and
$\alpha=2$ to compute $a_{0,2}^{bd}$ and $a_{2,2}^{bd}$:
\medbreak\qquad $a_{0,2}^{bd}(F,\DD)=\kappa_2(4\pi)^{-m/2}
\mathcal{I}^{bd}\{\Tr(-F_0\Id)\}$ \medbreak\qquad\qquad\qquad$
=\sqrt\pi(4\pi)^{-m/2} \mathcal{I}^{bd}\{\Tr(F_0\Id)\}$,
\medbreak\qquad
$a_{1,1}^{bd}(F,\DD)=\kappa_{0}(4\pi)^{-m/2}\mathcal{I}^{bd}\{\Tr(
-F_1\Id+\frac3{4}F_0L_{aa}\Id)\}$ \smallbreak\qquad\qquad\qquad$
=\frac{\sqrt\pi}8(4\pi)^{-m/2}\mathcal{I}^{bd}\{\Tr(
-4F_1\Id+3F_0L_{aa}\Id)\}$, \smallbreak\qquad
$a_{2,2}^{bd}(F,\DD)=\kappa_0(4\pi)^{-m/2}\mathcal{I}^{bd}\{\Tr(
-F_2\Id+\frac{\alpha-5}{2(\alpha-4)}F_1L_{aa}\Id
$\smallbreak\qquad\qquad\qquad$ +\frac{1}{6}F_0R_{amma}\Id
-\frac{\alpha-7}{8(\alpha-6)}F_0L_{aa}L_{bb}\Id
+\frac{\alpha-5}{4(\alpha-6)}F_0L_{ab}L_{ab}\Id$\smallbreak\qquad\qquad\qquad$
-\frac1{3(1-\alpha)}F_0R_{ijji}\Id-\frac2{1-\alpha}F_0E)\}|_{\alpha=2}$
\smallbreak\qquad\qquad\qquad$
=\sqrt\pi(4\pi)^{-m/2}\mathcal{I}^{bd}\{\Tr( -\frac{1}2F_2\Id
+\frac3{8}F_1L_{aa}\Id +\frac1{12}F_0R_{amma}\Id
$\smallbreak\qquad\qquad\qquad$ -\frac{5}{64}F_0L_{aa}L_{bb}\Id+
\frac3{32}F_0L_{ab}L_{ab}\Id+\frac16F_0R_{ijji}\Id+F_0E)\}$.
\medbreak
 We compute the remaining coefficients as follows. Let
$c(a_{i,\alpha}^{bd},A)$ be the coefficient of the monomial
$\mathcal{I}^{bd}\{\Tr(A)\}$ in
$(4\pi)^{m/2}a_{i,\alpha}^{bd}(F,\DD)$. It follows by Theorem
\ref{thm-1.5} that
$$
a_{0,1}^{bd}(F,\DD)={\textstyle\frac C2}\cdot(4\pi)^{-m/2}\mathcal{I}^{bd}\{\Tr(F_0\Id)\}\,.
$$
We may expand
$$
\textstyle\kappa_{\alpha-1}\frac{\alpha-4}{2(\alpha-3)}=\kappa_{\alpha-1}\frac{2(\alpha-3)+(2-\alpha)}{2(\alpha-3)}
=\textstyle\kappa_{\alpha-1}+\Gamma\left(\frac{4-\alpha}2\right)\frac1{2(\alpha-3)}\,.
$$
It now follows that
$$
\begin{array}{l}
\textstyle c(a_{1,2}^{bd},F_0L_{aa}\Id)=-{\textstyle\frac C2}-\frac12,\\
\textstyle a_{1,2}^{bd}(F,\DD)=(4\pi)^{-m/2}\mathcal{I}^{bd}\{\Tr({\textstyle\frac C2} F_1\Id-({\textstyle\frac
C2}+\frac12)F_0L_{aa}\Id)\}\,.
\vphantom{\vrule height 12pt}\end{array}$$
Many of the terms in $a_{2,1}^{bd}(F,\DD)$ are in fact regular at $\alpha=1$. We have:
$$\begin{array}{ll}
c(a_{2,1}^{bd},F_2\Id)=-\textstyle\frac12,&
c(a_{2,1}^{bd},F_1L_{aa}\Id)=\frac13,\\
c(a_{2,1}^{bd},R_{amma}\Id)=\textstyle\frac1{12},&
c(a_{2,1}^{bd},L_{aa}L_{bb}\Id)=-\textstyle\frac3{40},\vphantom{\vrule height 11pt}\\
c(a_{2,1}^{bd},L_{ab}L_{ab}\Id)=\textstyle\frac1{10}\vphantom{\vrule height 11pt}\,.
\end{array}$$
The terms involving $R_{ijji}$ and $E$ can be written in the form
$-\kappa_\alpha F_0(\frac16R_{ijji}\Id+E)$. Thus we may use the regularization of
$-F_0\Id$ in
$a_{0,1}^{bd}$ which was already computed to see:
$$
\begin{array}{l}
a_{2,1}^{bd}(F,\DD)=(4\pi)^{-m/2}\mathcal{I}^{bd}\left\{
\Tr\left(\textstyle-\frac12F_2\Id+\frac13F_1L_{aa}\Id+\frac1{12}R_{amma}\Id\right.\right.\\
\qquad\qquad\qquad\qquad\textstyle\left.\left.
-\frac3{40}L_{aa}L_{bb}\Id+\frac1{10}L_{ab}L_{ab}\Id+{\textstyle\frac C{12}} R_{ijji}\Id+{\textstyle\frac C2}
E\right)\right\}\,.
\vphantom{\vrule height 12pt}
\end{array}$$
This completes the derivation of Theorem \ref{thm-1.8} from Theorem \ref{thm-1.7}.\hfill\qedbox

\section{The pseudo-differential calculus}\label{sect-4}

We adopt the following notational conventions. Let $\al=(\alpha_1,...,\alpha_m)$ be a multi-index. We set:
$$\begin{array}{ll}
 |\al | = \alpha _1 + ...+\alpha _m,&
\al! = \alpha _1 ! \times ...\times \alpha_m !, \\
x^{\al} = x_1^{\alpha_1} \times ... \times x_m ^{\alpha_m}, &
d_x^\al = \left( \frac \partial {\partial x_1} \right)^{\alpha _1}
\times ... \times \left( \frac \partial {\partial x_m}
\right)^{\alpha _m}, \\ D_{\al} ^x = (-\sqrt{-1})^{|{\al}|} d_x^{\al}\,.
\end{array}$$
We apologize in advance for the slight notational confusion involved with using $\alpha$ to control the growth of $F$ and also
to using ${\al}$ as a multi-index.

We begin our discussion by reviewing the standard pseudo-differential computation of the resolvent on
a closed manifold without boundary and refer to \cite{ DW92, DOP82, G94,Gr68,Se69} for further details. Let $\DD$ be an
operator of Laplace type.  We want to construct the resolvent of
$\DD-\lambda$ for large $\lambda$ where we use Equation (\ref{eqn-1.b}) to express:
$$\DD=\sum_{|{\al}| \leq 2} a_{\al} (x) D_{\al} ^x.$$ For the symbol
$\sigma (\DD) (x,\xi)$ of $\DD$ this means
$$\sigma (\DD) (x,\xi ) = \sum_{|{\al}| \leq 2} a_{\al}(x) \xi ^{\al};$$
note that for the scalar Laplacian $\D$ the $0^{\operatorname{th}}$
term vanishes so  $a_{\vec0}=0$ in this setting. In the evaluation
of the heat equation asymptotics homogeneity properties of symbols
are relevant and it turns out that collecting terms according to
\begin{eqnarray} a_2 (x,\xi ,
\lambda ) &=& -
\lambda + \sum_{|\al | =2} a_\al (x) \xi ^\al , \nn\\
a_j (x,\xi ,\lambda ) &=& \sum_{|\al | = j} a_\al (x) \xi ^\al ,
\quad \quad j=0,1,\nn
\end{eqnarray}
 is fruitful. As a result, the symbol
$\sigma (\DD-\lambda)(x,\xi ,\lambda)$ can be written as
\begin{eqnarray}\sigma
(\DD-\lambda ) (x,\xi ,\lambda ) = \sum_{j=0}^2 a_j (x,\xi
,\lambda).\nn
\end{eqnarray}

 For the symbol of the resolvent of $\DD-\lambda$ we
make the Ansatz
\begin{eqnarray}\sigma ((\DD-\lambda )^{-1}) (x,\xi , \lambda )
\sim \sum_{l=0} ^\infty q_{-2-l} (x,\xi ,\lambda ) .\label{eqn-4.a}
\end{eqnarray}
In view of the formula for the symbol of a product we see that
$q_{-2-l}$ is determined algebraically by
\begin{eqnarray}
1 &=& a_2 (x,\xi , \lambda ) q_{-2} (x,\xi
,\lambda) , \label{eqn-4.b}\\
0 &=& \sum_{\stackrel{\al , j , l\leq k}{k=2+l+|\al |-j}}\frac 1
{\al !} [d_\xi^\al a_j (x,\xi , \lambda ) ] \,\, [ D_\al ^x
q_{-2-l} (x,\xi ,\lambda ) ] \quad \quad \mbox{for }k\geq
1.\label{eqn-4.c}
\end{eqnarray}
These symbols will play a crucial role in the proof of Theorem \ref{thm-1.6} in Section \ref{sect-5}.

We now specialize to the case where $\DD=\D$ is the scalar
Laplacian.
 For our present considerations we need
$q_{-2},q_{-3},q_{-4}$. For simplicity we skip the arguments in
the following summary of results -- we use in an essential fashion the fact that $\D$ is scalar.
We also use the convention that repeated indices are summed over.
We find (all greek indices will range over $\{1,2,...,m\}$):
\begin{eqnarray}
q_{-2} &=& a_2^{-1} , \nn\\q_{-3} &=& -a_2^{-1} \left[ a_1 q_{-2}
+
(D_\xi^\nu a_2 ) (\iD_\nu^x q_{-2})\right] , \nn\\
q_{-4} &=& - a_2 ^{-1} \left[ a_0 q_{-2} + a_1 q_{-3} +(D_\xi^\nu
a_1 ) (\iD_\nu^x q_{-2} ) \right.\nn\\
& &\left. +(D_\xi ^\nu a_2 ) (\iD_\nu^x q_{-3} ) - {\textstyle\frac12}
(D_\xi^{\nu\mu} a_2 ) (D_{\nu\mu} ^x q_{-2} ) \right]. \nn
\end{eqnarray}
 For
later use it will be advantageous to express the results in terms
of $q_{-2}^n$. We then have
\begin{eqnarray} q_{-3} &=& -a_1 q_{-2}^2 +
c_{-3,3} q_{-2}^3, \nn\\
q_{-4} &=& -a_0 q_{-2}^2 + c_{-4,3} q_{-2}^3 + c_{-4,4} q_{-2}^4 +
c_{-4,5} q_{-2}^5 , \nn
\end{eqnarray}
 where
\begin{eqnarray} c_{-3,3} &=& -\sqrt{-1}\,(\partial
_\xi ^\nu a_2) (\partial
_\nu^x a_2) ,\nn\\
c_{-4,3} &=& a_1^2 - \sqrt{-1}\,(\partial _\xi ^\nu a_1) (\partial _\nu ^x
a_2 ) - \sqrt{-1}\,(\partial_\xi ^\nu a_2) (\partial _\nu ^x a_1 ) - {\textstyle\frac12} (\partial _\xi ^{\nu\mu} a_2 )
(\partial _{\nu\mu} ^x a_2 ) ,
\nn\\
c_{-4,4} &=& -3 a_1 c_{-3,3} + \sqrt{-1}\,(\partial _\xi^\nu a_2 )
(\partial _\nu ^x c_{-3,3} ) + (\partial _\xi ^{\nu\mu} a_2 )
(\partial _\nu ^x a_2 ) (\partial _\mu^x a_2) , \nn\\
c_{-4,5} &=& 3 c_{-3,3}^2.\nn
\end{eqnarray}
 The relevant operator for our
considerations is
$$ \D-\lambda = - g^{\mu\nu} \partial _\mu
\partial _\nu + b^\mu \partial _\mu - \lambda
$$
where we have changed notation slightly from that used previously. For the symbols this gives
\begin{eqnarray*} a_2 (x,\xi , \lambda ) &=&
g^{\mu\nu} \xi_\mu \xi _\nu - \lambda \equiv |\xi|^2 -\lambda ,
\nn\\
a_1 (x,\xi ,\lambda ) &=& \sqrt{-1} b^\mu \xi_\mu , \quad\text{and}\quad
a_0 (x,\xi ,\lambda )=0\,.
\end{eqnarray*}
 To state results for $q_{-2}$,
$q_{-3}$ and $q_{-4}$ for this operator $D$ we will as usual raise
and lower indices using the inverse metric and the metric.
Furthermore, ``$,$" denotes partial differentiation. One computes easily that:
\begin{eqnarray*}
&&q_{-2} (x ,\xi ,\lambda )= \textstyle\frac 1 {|\xi|^2 -\lambda } ,
\\
&&q_{-3} (x , \xi , \lambda) =\textstyle - \frac {1} { (|\xi|^2
-\lambda )^2}\sqrt{-1}b^\mu \xi_\mu  - \frac {1}{(|\xi|^2 -\lambda )^3} 2\sqrt{-1}g^{\sigma  \gamma }_{,\nu} \xi^\nu \xi_\sigma
\xi_\gamma,\\
&&q_{-4} (x, \xi , \lambda )=
\\
& &\quad\textstyle\frac 1 {(|\xi|^2 -\lambda )^3} \left\{ - b^\mu b^\nu \xi_\mu
\xi _\nu + b^{\nu } g^{\sigma   \beta}_{,\nu} \xi _\sigma  \xi _\beta + 2
b^\sigma  _{,\nu} g^{\nu\beta} \xi_\beta \xi_\sigma  - g^{\sigma  \beta}
_{,\nu\mu} g^{\nu\mu} \xi _\sigma  \xi _\beta \right\} \\
& &\quad\textstyle+\frac 1 {(|\xi|^2 -\lambda )^4} \left\{ -6 b^\mu g^{\sigma  \gamma}
_{,\nu} g^{\nu\beta} \xi _\mu \xi _\beta  \xi_\sigma  \xi _\gamma + 4
g^{\sigma  \gamma }_{,\beta \nu} g^{\beta \mu} g^{\nu\delta} \xi _\mu
\xi _\sigma  \xi _\gamma \xi _\delta \right.\\
& &\qquad\left.+4 g^{\sigma  \gamma }_{,\beta } g^{\beta
\mu}_{,\nu} g^{\nu\delta} \xi _\mu \xi _\sigma  \xi _\gamma \xi
_\delta + 2 g^{\sigma  \beta} _{,\nu} g^{\gamma \delta } _{,\mu}
g^{\nu\mu} \xi
_\sigma  \xi _\beta \xi _\delta \xi _\gamma\right\}\nn\\
& &\quad\textstyle+\frac 1 {(|\xi|^2 -\lambda )^5} \left\{ -12 g^{\sigma  \gamma}
_{,\nu} g^{\nu\beta} g^{\delta\tau}_{,\mu} g^{\mu\rho} \xi_\beta
\xi_\sigma  \xi_\gamma\xi_\rho \xi_\delta \xi_\tau\right\}.\nn
\end{eqnarray*}

If the manifold has a boundary the expansion (\ref{eqn-4.a}) has to be
augmented by a boundary correction. To formulate the conditions to
be satisfied by the boundary correction we expand about $r=0$. We adopt the notation
established in Section \ref{sect-1.7} and expand
$$ds^2_M=g_{\sigma\varrho}(y,r)dy^\sigma\circ dy^\varrho+dr^2\quad\text{on}\quad\mathcal{C}_\varepsilon\,.$$
 The coordinate $y$ locally parameterizes the boundary, and
$r$ is the geodesic distance to the boundary, so $x=(y,r)$. A
tilde above any quantity will indicate that it is to be evaluated
at the boundary, that is at $r=0$. Furthermore, we use $\xi
=(\omega ,\tau )$.

We find $$\D-\lambda = \sum_{k=0}^\infty \frac 1 {k!} r^k
\sum_{|\al | \leq 2} \frac{\partial ^k} {\partial r^k} a_\al
(y,r)\bigg |_{r=0} D_{y,r}^\al $$ with the notation $$D_{y,r}^\al
= \left(\prod_{i=1}^{m-1} D_{y_i}^{\alpha_i}\right) D_r^{\alpha
_m}.$$ Introducing
\begin{eqnarray} a_j ( \vard ) =\left\{\begin{array}{ll}
\sum_{|\al | =j} a_\al (y,r) \left( \prod_{i=1}^{m-1}
\omega_i^{\alpha _i} \right) D_r^{\alpha_m}&\text{for}\quad j=0,1,\\
\sum_{|\al | =2} a_\al (y,r) \left( \prod_{i=1}^{m-1}
\omega_i^{\alpha _i} \right) D_r^{\alpha_m}-\lambda&\text{for}\quad j=2\end{array}\right. \nn
\end{eqnarray}
 we define the
partial symbol
\begin{eqnarray} \sigma ' (\D-\lambda ) = \sum_{k=0}^\infty
\frac 1 {k!} r^k \sum_{j=0}^2 \frac{\partial ^k } {\partial r^k}
a_j (\vard ) \bigg|_{r=0} .\nn
\end{eqnarray}
 As it turns out, the symbols
\begin{eqnarray}
a^{(j)} ( \vard ) = \sum_{l=0}^2
\sum_{\stackrel{k=0}{l-k=j}}^\infty \frac 1 {k!} r^k
\frac{\partial ^k}{\partial r^k} a_l ( \vard )\bigg |_{r=0} \nn
\end{eqnarray}
have suitable homogeneity properties and using these symbols we
write
\begin{eqnarray} \sigma ' (\D-\lambda ) = \sum_{j=-\infty} ^2 a^{(j)}
(\vard ) . \nn
\end{eqnarray}
 We write the symbol of the resolvent as
\begin{equation}\label{eqn-4.d}
\begin{array}{l}
\displaystyle\sigma
((\D-\lambda )^{-1}) (\var ) \\
\displaystyle=\sum_{j=0}^\infty q_{-2-j} (y,r,\omega ,\tau
,\lambda ) - e^{-\sqrt{-1}\tau r}\sum_{j=0}^\infty h_{-2-j} (\var ) ,
\end{array}\end{equation}
 where the second term is the boundary
correction. The factor $e^{-\sqrt{-1}\tau r}$ appears because the operator constructed from
these terms is the $Op^\prime(h)$ in \cite{Se69}, and $Op^\prime(h)=Op(he^{-\sqrt{-1}\tau r})$. This shows
\begin{eqnarray} \sigma ' (\D-\lambda )\circ
\sum_{j=0}^\infty h_{-2-j} (\var ) =0 .\nn
\end{eqnarray}
Here $\circ$ denotes the symbol product on $\mathbb{R}^{m-1}$. Analogously to
Equations (\ref{eqn-4.b}) and (\ref{eqn-4.c}) this equation leads to the
differential equations
\begin{eqnarray} 0 &=& a^{(2)} ( \vard ) h_{-2} (\var ),
\nn\\ 0 &=& a^{(2)} (\vard ) h_{-2-j} (\var ) \nn\\
& &+ \sum_{\stackrel{\al , k , l<j}{j=l+2+|\al | -k}} \frac
1 {\al!} \left[ D_\omega ^\al a^{(k)} ( \vard ) \right]
\left[ (\iD^y)_\al h_{-2-l} (\var ) \right] .\nn
\end{eqnarray}
 For the
present considerations we need $h_{-2-j}$ for $j=0,1,2$, and we
have more
explicitly (repeated letters $a,b,c,...$ run over tangential coordinates $\{1,2,...,m-1\}$)
\begin{eqnarray*} 0 &=& a^{(2)} (\vard ) h_{-2} (\var ), \nn\\
0 &=&a^{(2)} (\vard ) h_{-3} (\var ) + a^{(1)} (\vard ) h_{-2}
(\var ) \nn\\
& & + \left[ D_\omega ^b a^{(2)} (\vard ) \right] \left[ (\iD^y)_b
h_{-2} (\var ) \right], \nn\\
0&=&a^{(2)} (\vard ) h_{-4} (\var ) + a^{(0)} (\vard ) h_{-2}
(\var ) \nn\\
& &+ \left[ D_\omega ^b a^{(1)} (\vard ) \right] \left[ (\iD^y)_b
h_{-2} (\var ) \right] \nn\\
& &+ {\textstyle\frac12} \left[ D_\omega ^{bc} a^{(2)}
(\vard ) \right] \left[ (\iD^y )_{bc} h_{-2} (\var ) \right] \nn\\
& &+ a^{(1)} (\vard ) h_{-3} (\var ) \\&&+ \left[ D_\omega ^b a^{(2)}
(\vard ) \right] \left[ (\iD^y ) _b h_{-3} (\var ) \right] .\nn
\end{eqnarray*}
The relevant equations for $a^{(i)} (\vard )$, $i=0,1,2$ are
\begin{eqnarray*}
a^{(2)} (\vard ) &=& a_2 (\vard ) |_{r=0} \nn\\
&=&\gt \omega_a \omega_b + D_r^2 -\lambda , \nn\\ a^{(1)} (\vard )
&=& r (\partial_r a_2 (\vard ))|_{r=0} + a_1 (\vard ) |_{r=0}
\nn\\
&=& r \gtr \omega_a \omega_b + \sqrt{-1} \tilde b^a \omega _a + \sqrt{-1} \tilde
b^r D_r , \nn\\
a^{(0)} (\vard ) &=& {\textstyle\frac12} r^2 (\partial_r^2 a_2 (\vard )
)|_{r=0} + r (\partial _r a_1 (\vard ) )|_{r=0} \\&&+ a_0 (\vard )
|_{r=0} \nn\\
&=& {\textstyle\frac12} r^2 \gtrr \omega_a \omega_b + r\sqrt{-1}\,\tilde b^a_{,r}
\omega_a + r \sqrt{-1}\, \tilde b^r _{,r} D_r + \tilde c .\nn
\end{eqnarray*}

The differential equations have to be augmented by a growth
condition
\begin{eqnarray} h_{-2-j} (\var )
\to 0 \quad \quad \mbox{as}\quad r\to\infty \label{eqn-4.e},
\end{eqnarray}
and an initial condition corresponding to the Dirichlet boundary
condition
\begin{equation}\label{eqn-4.f}
h_{-2-j} (\var ) |_{r=0} = q_{-2-j} (\var ) |_{r=0}\,.
\end{equation}
 Once the
symbols $h_{-2-j}$ have been determined, their contribution to the
asymptotics of the trace of the heat kernel follows from multiple
integration. As before, we suppose $r^\alpha F\in C^\infty(\mathcal{C}_{\varepsilon})$.
The contribution reads $$\sum_{l=0}^\infty t^{\frac {1-\alpha -m}
2 } t^{\frac l 2} \int_{\partial M} \eta _{\frac l 2}
(y,F,\D)  dy$$ with
\begin{eqnarray}&& \eta _{\frac l 2} (y,F,\D )= \frac 1
{(2\pi)^{m+1}} \sum_{j+k=l}\nn
\int_{\mathbb{R}^{m-1}} d\omega
\int_{-\infty }^\infty ds\nn\\&&\quad\times \int_0^\infty d\bar r
e^{\sqrt{-1}\,s} \left( - \int_\gamma d\tau e^{-\sqrt{-1}\,\tau \bar r}\right)
h_{-2-j} (y,\bar r,\omega,\tau, -\sqrt{-1}\,s) \bar r ^{k-\alpha}F_k(y),\label {eqn-4.g}
\end{eqnarray}
 where $\gamma$ is anticlockwise enclosing
the poles of $h_{-2-j}$ in the lower half-plane. The integral with
respect to $s$ is the contour integral transforming the resolvent to
the heat kernel; see Section 5. Note that from (\ref{eqn-4.d}) the
contribution to the heat kernel is {\bf minus} the above.

As will become clear in the following, with $\Lambda =
\sqrt{|\omega|^2 +\sqrt{-1}\,s}$, we need integrals of the type
\medbreak\noindent\centerline{$\displaystyle
T_{ab...}^{kljn} \equiv \int_{\mathbb{R}^{m-1}} d\omega
\int_{-\infty }^\infty ds \int_0^\infty d\bar r
e^{\sqrt{-1}\,s} \left( - \int_\gamma d\tau e^{-\sqrt{-1}\,\tau \bar r}\right)
\frac{\tau ^k \bar r ^{l-\alpha} \omega_a \omega_b
...}{\Lambda^j (\tau^2 + \Lambda^2 ) ^n} e^{-\bar r\Lambda }
.\nn$
}\medbreak\noindent
 The $\tau$ integration can be done using
\begin{eqnarray}\int_\gamma d\tau e^{-\sqrt{-1}\,\tau \bar r} \frac{\tau
^k}{(\tau ^2 + \Lambda^2)^l} &=& \frac{(\sqrt{-1})^k (-1)^{l+k} \pi
}{(l-1)!} \left( \frac 1 {2\Lambda} \frac d
{d\Lambda}\right)^{l-1} \left[ \Lambda^{k-1} e^{-\bar r
\Lambda}\right] .\nn
\end{eqnarray}
 So
\begin{eqnarray*} T_{ab...}^{kljn} &=& {\textstyle\frac{(\sqrt{-1})^k
(-1) ^{n+k+1} \pi}{(n-1)!}}  \int_{\mathbb{R}^{m-1}} d\omega
\int_{-\infty}^\infty ds\\&&\quad\times \int_0^\infty d\bar r
e^{\sqrt{-1}\,s} \bar r ^{l-\alpha} {\textstyle\frac{\omega_a \omega_b ...} {\Lambda^j}}
e^{-\bar r \Lambda} \left( {\textstyle\frac 1 {2\Lambda} \frac d {d\Lambda}}
\right)^{n-1} \left[ \Lambda^{k-1} e^{-\bar r
\Lambda}\right].\nn
\end{eqnarray*}
 Performing the $\Lambda$-differentiation,
different $\bar r$-dependent functions would occur. It is
therefore desirable to first perform the $\bar r$-integration
before performing the $\Lambda$-derivatives explicitly. This is
achieved by noting that ($z=\Lambda$ has to be put after
the $\Lambda$ differentiation has been performed)
\begin{eqnarray*} T_{ab...}^{kljn} &=& \left.\frac{(\sqrt{-1})^k (-1)
^{n+k+1} \pi}{(n-1)!} \int_{\mathbb{R}^{m-1}} d\omega
\int_{-\infty}^\infty ds e^{\sqrt{-1}\,s} \frac{\omega_a \omega_b
...} {\Lambda^j}\right.\\&&\quad\left.\times\left( \frac 1 {2\Lambda} \frac d {d\Lambda}
\right)^{n-1} \Lambda^{k-1} \int_0^\infty d\bar r \bar
r^{l-\alpha} e^{-\bar r (\Lambda + z)} \right|_{z=\Lambda}\nn\\
&=&\left.\frac{(\sqrt{-1})^k (-1) ^{n+k+1} \pi}{(n-1)!}\Gamma(l+1-\alpha )
\int_{\mathbb{R}^{m-1}} d\omega \int_{-\infty}^\infty ds
e^{\sqrt{-1}\,s} \frac{\omega_a \omega_b ...} {\Lambda^j}\right.\\&&\qquad\left. \left( \frac 1
{2\Lambda} \frac d {d\Lambda} \right)^{n-1}
\frac{\Lambda^{k-1}}{(\Lambda +
z)^{l+1-\alpha}}\right|_{z=\Lambda}.\nn
\end{eqnarray*}
 We can proceed in general by introducing numerical multipliers $c_{nkl}$
according to
\begin{eqnarray} \left. \left( \frac 1 {2\Lambda} \frac d
{d\Lambda}\right)^{n-1} \frac{\Lambda^{k-1}}{(\Lambda
+z)^{l+1-\alpha }} \right|_{z=\Lambda} = c_{nkl} \frac 1
{\Lambda^{l+2n-k-\alpha}}.\nn
\end{eqnarray}
 The $s$-integration is then
performed using $$\int_{-\infty}^\infty ds
\frac{e^{\sqrt{-1}\,s}}{(|\omega|^2 +\sqrt{-1}\,s )^\beta} = \frac{2\pi} {\Gamma (\beta
)} e^{-|\omega|^2} .$$ The final $\omega$-integrations follow from
\begin{eqnarray} C(y) & \equiv & \int_{\mathbb{R} ^{m-1}} d\omega e^{-\gt
\omega_a \omega_b + \sqrt{-1}\, y^a \omega_a} = \pi^{\frac{m-1} 2}
\sqrt{\tilde g} e^{-\frac {\tilde g _{ab} y^a y^b} 4} , \nn
\end{eqnarray}
 by
observing that
\begin{eqnarray} \int_{\mathbb{R} ^{m-1}} d\omega \,\,
\omega_{a_1} \omega _{a_2} ... \omega_{a_r} e^{-\gt \omega_a
\omega_b} &=& \left. \left( \frac 1 {\sqrt{-1}} \right)^r \frac \partial
{\partial y^{a_1}} \cdot\cdot\cdot \frac \partial {\partial
y^{a_r}} C(y) \right|_{y=0}.\nn
\end{eqnarray}
 In particular
\begin{eqnarray}
\int_{\mathbb{R} ^{m-1}} d\omega \,\,
e^{-|\omega|^2} &=& \pi^{\frac{m-1} 2} \sqrt {\tilde g}, \nn\\
\int_{\mathbb{R} ^{m-1}} d\omega \,\,\omega_a \omega_b
e^{-|\omega|^2} &=& {\textstyle\frac12} \pi^{\frac{m-1} 2} \sqrt{\tilde g}
\tilde g_{ab}\nn,\\
\int_{\mathbb{R} ^{m-1}} d\omega \,\,\omega_a
\omega_b\omega_c\omega_d e^{-|\omega|^2} &=& {\textstyle\frac14}
\pi^{\frac{m-1} 2} \sqrt{\tilde g} \left( \tilde g_{ab} \tilde
g_{cd} + \tilde g _{ac} \tilde g _{bd} + \tilde g_{ad} \tilde
g_{bc} \right).\nn
\end{eqnarray}
 Introducing the numerical multipliers
$d_{kljn}$ according to
\begin{eqnarray} d_{kljn} = \frac {2 (\sqrt{-1})^k (-1)^{n+k+1}
\pi^2 \Gamma (l+1-\alpha ) c_{nkl}}{(n-1)! \Gamma \left( \frac{
j+l-k-\alpha } 2 +n \right)} , \nn
\end{eqnarray}
 we obtain the
compact-looking answers
\begin{eqnarray} T_{ab...}^{kljn} = d_{kljn}
\int_{\mathbb{R} ^{m-1}} d\omega \,\, \omega_a \omega_b ...
e^{-|\omega|^2},\nn
\end{eqnarray}
 where the last $\omega$-integration is
performed with the above results.

Note that the numerical multipliers $d_{kljn}$ are easily
determined using an algebraic computer program. Therefore, all
appearing integrals can be very easily obtained.

Let us apply this formalism explicitly to the leading orders, and
we start with $h_{-2} (\var ).$ The relevant differential equation
reads $$(\partial_r^2 - \Lambda ^2) h_{-2} (\var ) =0,$$ which has
the general solution $$h_{-2} (\var ) = A e^{-r \Lambda} + B e^{r
\Lambda} .$$ The asymptotic condition (\ref{eqn-4.e}) on the symbol
as $r\to\infty$ imposes $B=0$. The initial condition $h_{-2}
|_{r=0} = q_{-2} |_{r=0}$ gives $$A = \frac 1 {\tau^2 +
\Lambda^2}.$$ Putting the information together we have obtained
$$h_{-2} (\var ) = \frac 1 {\tau ^2 + \Lambda^2} e^{-r \Lambda}
.$$ Performing the relevant integrals, with the notation $$\int dI
= \int_{\mathbb{R} ^{m-1}} d\omega \int_{-\infty}
^\infty ds \int_0^\infty d\bar r e^{\sqrt{-1}\,s}\left( -
\int_\gamma d\tau e^{-\sqrt{-1}\,\tau \bar r}\right) \bar r
^{-\alpha} ,$$ produces
\begin{eqnarray*}
&&\int dI h_{-2} ( y,\bar r,\omega,\tau ,-\sqrt{-1}\,s ) = d_{0001} \pi^{\frac{m-1} 2 } \sqrt
{\tilde g} = \frac{2^\alpha \pi^2 \Gamma (1-\alpha) } { \Gamma
\left( 1-\frac \alpha 2 \right)} \pi ^{\frac{m-1} 2 } \sqrt
{\tilde g}\\&=& \pi \Gamma \left( \frac{1-\alpha} 2 \right) \pi^{m/2}
\sqrt {\tilde g}\,.
\end{eqnarray*} Taking into account the prefactor in
(\ref {eqn-4.g}), this confirms the value of $\bar\kappa_\alpha$ in Lemma \ref{lem-1.9}.

In the next order we obtain
\begin{eqnarray} (\partial _r^2 - \Lambda^2)
h_{-3} (\var ) = (E+U_1) e^{-r \lambda} + (F+U_2) r e^{-r \Lambda}
, \nn
\end{eqnarray}
 where
\begin{eqnarray*} E&=& -\frac{ \tilde b_r \Lambda}{\tl } ,
\quad \quad F= \frac{\gtr \omega_a \omega_b} {\tl} , \nn\\
U_1(\omega ) &=& \frac{ \sqrt{-1}\, \tilde b^a \omega_a} {\tl }  + \frac {2\sqrt{-1}\,
\tilde g^{ac} _{,b} \omega^b\omega_a \omega_c} {(\tl )^2}, \quad
U_2 (\omega ) = \frac{ \sqrt{-1}\, \tilde g^{ac} _{,b} \omega^b \omega_a
\omega_c}{(\tl ) \Lambda}  .
\end{eqnarray*}
 Note, for later
arguments, that $U_1(\omega )$ and $U_2 (\omega )$ are odd functions
in $\omega$. Furthermore, for the scalar Laplacian at hand $b^a =
g^{bc} \Gamma_{bc}{}^a$; thus they contain only {\it tangential}
derivatives of the metric.

Using for example the annihilator method, we write down the
general form of the solution to this differential equation as
$$h_{-3} (\var ) = c_1 e^{-r \Lambda} + c_2 r e^{-r \Lambda} + c_3
r^2 e^{-r \Lambda} + c_4 e^{r\Lambda}.$$ From the asymptotic
condition (\ref{eqn-4.e}) we conclude $c_4 =0$. From the initial
condition given in Equation (\ref{eqn-4.f}) we obtain
\begin{eqnarray} c_1 = - \frac{\sqrt{-1}\,\tilde b^a \omega} {(\tl
)^2} - \frac {\sqrt{-1}\,\tilde b^r \tau}{(\tl )^2} - \frac{ 2\sqrt{-1}\, \tilde
g^{ab} _{,c} \omega^c \omega_a \omega_b}{(\tl )^3} - \frac {2\sqrt{-1}\,\gtr
\tau \omega_a \omega_b} {(\tl )^3} .\nn
\end{eqnarray}
 From the differential
equation we derive
\begin{eqnarray} c_2 &=& - {\textstyle\frac 1 {4 \Lambda^2}} (F+U_2) -
{\textstyle\frac 1 {2\Lambda}} ( E+U_1) , \nn\\
c_3 &=& - {\textstyle\frac 1 {4 \Lambda}} (F+U_2) .\nn
\end{eqnarray}
 Collecting the
available information, we see
\begin{eqnarray} h_{-3} (\var ) = D e^{-r
\Lambda} + Br e^{-r \Lambda} + C r^2 e^{-r \Lambda} + O(\omega ) ,
\nn
\end{eqnarray}
 with
\begin{eqnarray} D&=& - \frac{ \sqrt{-1}\, \tilde b^r \tau }{(\tl )^2} -
\frac{ 2\sqrt{-1}\, \gtr \tau \omega_a \omega_b} {(\tl )^3} , \nn\\
B &=& - \frac{ \gtr \omega_a \omega_b}{4 \Lambda^2 (\tl ) } +
\frac{ \tilde b^r } {2 (\tl ) } , \nn\\
C&=& - \frac{ \gtr \omega_a \omega_b} { 4 \Lambda (\tl )} ,
\nn
\end{eqnarray}
 and where $O(\omega )$ is an odd function in $\omega$.
Furthermore, $O(\omega )$ contains only tangential derivatives of
the metric. We next perform the multiple integrals; note, odd
functions in $\omega$ do not contribute. We obtain
\begin{eqnarray*} &&\int dI
h_{-3} (y,\bar r ,\omega,\tau,-\sqrt{-1}\,s )\\ &=& \pi^{\frac{m-1} 2 } \sqrt
{\tilde g} \gtr \tilde g_{ab} \left\{ -\textstyle \frac {\sqrt{-1}} 2 d_{1002} + {\textstyle\frac14} d_{0101} -
\sqrt{-1}\, d_{1003} - {\textstyle\frac18} d_{0121} - {\textstyle\frac18} d_{0211} \right\}\nn\\
&=&\frac{\pi (\alpha -4)}{4 (3-\alpha)} \Gamma \left(
\frac{2-\alpha} 2 \right) \pi^{m/2} \sqrt {\tilde g} \gtr \tilde
g_{ab} .\nn
\end{eqnarray*}
This confirms the value of $\kappa^1_\alpha$ in Lemma \ref{lem-1.9} after taking into account the prefactor in (\ref
{eqn-4.g}) and the fact that $\gtr \tilde g_{ab}= - \tilde g^{ab} \tgr = 2
g^{ab}L_{ab}$.

Up to this point the calculation can be considered a warm up for
the next order. We would like to
determine the universal coefficients of the geometric invariants
$R_{amma}$, $L_{aa} L_{bb}$ and $L_{ab} L_{ab}$. In terms of the
metric the last two are determined by
$$L_{ab} = - {\textstyle\frac12} \tgr\,.$$
Using the Christoffel symbols
$$\Gamma_{jk}{} ^i = {\textstyle\frac12}
g^{il} \left( g_{lj,k} + g_{kl,j} - g_{jk,l} \right),$$
and taking into account that with our sign convention the scalar curvature is given by the contraction $g^{jk}R_{ijk}{}^i$,
we may expand the Riemann curvature tensor in the form:
$$
R_{ijk}{}^l=\Gamma_{jk}{}^l{}_{,i}-\Gamma_{ik}{}^l{}_{,j}+\Gamma_{in}{}^l\Gamma_{jk}{}^n-\Gamma_{jn}{}^l\Gamma_{ik}{}^n\,.
$$
The normal projection of the Riemann curvature tensor
reads
\begin{eqnarray} \tilde R_{amma} &=& - {\textstyle\frac12} \tilde g^{ac} _{,r}
\tilde g_{ac,r} - {\textstyle\frac12} \tilde g^{ac}\tilde  g_{ac,rr} - {\textstyle\frac14} \tilde g^{bc}_{,r} \tilde g^{ad} _{,r} \tilde g_{ca}
\tilde g_{bd} \nn\\
&=& {\textstyle\frac14} \tilde g^{ab} \tilde g^{cd} \tilde g_{ac,r} \tilde
g_{bd,r} - {\textstyle\frac12} \tilde g^{ac}\tilde  g_{ac,rr}.\nn
\end{eqnarray}
 The
above results suggest a strategy for the calculation. It suffices to
consider the special case where the metric is independent of $y$. As
a consequence, our answer will have the form
\begin{equation}\label{eqn-4.h}(4 \pi )^{-m/2}
\left\{ A \tilde g^{ac} \tilde g_{ac,rr} + B\tilde  g^{ab} \tilde
g^{cd}\tilde  g_{ac,r}\tilde g_{bd,r} + C \tilde g^{ab} \tilde
g^{cd} \tilde g_{ab,r}\tilde g_{cd,r} \right\}
.
\end{equation}
 This has to be compared with the terms in
$a_{2,\alpha}^{bd} (F,\Delta_M)$ that possibly contribute to these geometric
invariants. In detail one can show these terms are (mod terms with
tangential derivatives of the metric)
\begin{eqnarray} -{\textstyle\frac1{6}}
\kappa_\alpha \tilde R+ \kappa_\alpha^3\tilde  R_{amma} +
\kappa_\alpha^4 L_{aa} L_{bb} + \kappa_\alpha^5 L_{ab} L_{ab} &=&
\nn\\
& &\hspace{-6.0cm}\tilde g^{ac} \tilde g_{ac,rr} \left( {\textstyle\frac1{6}}
\kappa_\alpha - {\textstyle\frac12} \kappa_\alpha^3 \right) + \tilde g^{ab}
\tilde g^{cd} \tilde g_{ab,r}\tilde g_{cd,r} \left( {\textstyle\frac1{24}}
\kappa_\alpha + {\textstyle\frac14}
\kappa_\alpha^4 \right) \nn\\
& &\hspace{-6.0cm} + \tilde g^{ab} \tilde g^{cd}\tilde g_{ac,r}
\tilde g_{bd,r} \left( - {\textstyle\frac18} \kappa_\alpha + {\textstyle\frac14}
\kappa_\alpha ^3 +{\textstyle\frac14} \kappa_\alpha^5\right) .\nn
\end{eqnarray}
 So
once we know $A,B,C$, we can deduce
\begin{equation}\label{eqn-4.i}
\kappa_\alpha^3 = -2
\left( A - {\textstyle\frac1{6}} \kappa_\alpha \right),\quad \kappa_\alpha^4
= 4 \left( C - {\textstyle\frac1{24}} \kappa_\alpha
\right),\quad
\kappa_\alpha^5 = 4 \left(B+{\textstyle\frac18} \kappa_\alpha - {\textstyle\frac14}
\kappa_\alpha ^3 \right) \,.
\end{equation}
 In summary, when writing down
the differential equation for $h_{-4} (\var )$, we can neglect all
terms that are odd in $\omega$ as well as all terms that contain
tangential derivatives of the metric. We obtain (up to irrelevant
terms)
\begin{eqnarray} (\partial_r^2 -\Lambda ^2) h_{-4} (\var ) = A e^{-r
\Lambda} + B r e^{-r \Lambda } + C r^2 e^{-r \Lambda} + D r^3
e^{-r \Lambda } ,\nn
\end{eqnarray}
 where
\begin{eqnarray*} A &=& \textstyle\frac{ \sqrt{-1}\, \tilde b^r
\tilde b^r  \Lambda\tau }{(\tl )^2} + \frac{2\sqrt{-1}\,\tilde b^r \gtr
\Lambda \tau \omega_a \omega_b}{(\tl )^3 } - \frac {\tilde b^r
\gtr \omega_a \omega_b}{4 \Lambda^2 (\tl )} + \frac{ \tilde b^r
\tilde b^r}{2 (\tl )} , \nn\\
B &=& \textstyle- \frac{ \tilde b^r _{,r} \Lambda} {\tl } - \frac{ \sqrt{-1}\, \gtr
\tilde b^r \tau \omega_a \omega_b} {(\tl )^2} - \frac {2\sqrt{-1}\, \gtr
\tilde g^{cd} _{,r} \tau \omega_a \omega_b \omega_c \omega_d}{(\tl
)^3} - \frac{ \gtr \tilde b^r \omega_a \omega_b}{4 \Lambda (\tl )}
- \frac{ \tilde b^r \tilde b^r \Lambda }{2 (\tl )} ,\nn\\
C &=& \textstyle\frac {\gtrr \omega_a \omega_b}{2(\tl ) } - \frac{\gtr
\tilde g^{cd} _{,r} \omega_a \omega_b \omega_c \omega_d}{4 \Lambda
^2 (\tl )} +\frac{3\gtr \tilde b^r \omega_a \omega_b} {4 (\tl )} , \nn\\
D &=& \textstyle- \frac{ \gtr \tilde g^{cd} _{,r} \omega_a \omega_b \omega_c
\omega_d } {4 \Lambda (\tl )} .\nn
\end{eqnarray*}
 So the solution has the
form, taking into account the asymptotic behavior (\ref{eqn-4.e}),
\begin{eqnarray} h_{-4} (\var ) &=& \tilde \alpha e^{-r \Lambda} + \beta r
e^{-r \Lambda} + \gamma r^2 e^{-r \Lambda} + \delta r^3 e^{-r
\Lambda} + \epsilon r^4 e^{-r \Lambda} .\nn
\end{eqnarray}
 From the initial
condition $\tilde \alpha = q_{-4} (\var ) |_{r=0}$ we obtain, up
to irrelevant terms,
\begin{eqnarray} \tilde \alpha &=&\textstyle \frac 1 {(\tl )^3}
\left\{ - \tilde b^r \tilde b^r \tau ^2 + \tilde b^r \gtr \omega_a
\omega_b + 2
\tilde b^r_{,r} \tau ^2 - \gtrr \omega_a \omega_b\right\} \nn\\
&+&\textstyle \frac 1 {(\tl )^4 } \left\{ - 6 \tilde b^r \gtr \tau^2
\omega_a \omega_b + 4 \gtrr \tau^2 \omega_a \omega_b + 2 \gtr
\tilde g^{cd} _{,r} \omega_a \omega_b \omega _c \omega_d\right\}
\nn\\
&+&\textstyle \frac 1 {(\tl )^5 } \left\{ -12 \gtr \tilde g^{cd}_{,r}
\omega_a \omega_b \omega_c\omega_d \tau^2\right\} .\nn
\end{eqnarray}
 From
the differential equation we obtain the conditions
$$\begin{array}{ll} A = - 2
\Lambda \beta + 2 \gamma,&
B = - 4 \Lambda \gamma + 6 \delta ,\\
C = - 6 \Lambda \delta + 12 \epsilon,&
D = - 8 \epsilon \Lambda \,.\end{array}$$
 This determines the numerical
multipliers $\beta$, $\gamma$, $\delta$ and $\epsilon$ to be
\begin{eqnarray}
\beta &=& \textstyle- {\textstyle\frac38} \frac D {\Lambda^4} - {\textstyle\frac14} \frac C
{\Lambda^3} - {\textstyle\frac14} \frac B {\Lambda^2} - {\textstyle\frac12} \frac A
\Lambda , \nn\\
\gamma &=& \textstyle- {\textstyle\frac38} \frac D {\Lambda^3} - {\textstyle\frac14} \frac C
{\Lambda^2} - {\textstyle\frac14} \frac B \Lambda , \nn\\
\delta &=& \textstyle- {\textstyle\frac14} \frac D {\Lambda^2} - {\textstyle\frac1{6}} \frac C
\Lambda , \nn\\
\epsilon &=&\textstyle - {\textstyle\frac18} \frac D \Lambda .\nn
\end{eqnarray}
 For the
Laplacian on the manifold $M$ we have
\begin{eqnarray}
 \tilde b ^r &=& - {\textstyle\frac12} \tilde g^{ab} \tgr ,\nn\\
 \tilde b^r \tilde b^r &=& {\textstyle\frac14} \tilde g^{ab} \tilde g^{cd}
 \tgr \tilde g_{cd,r} ,\nn\\
 \tilde b^r \tilde g_{ab} \gtr &=& {\textstyle\frac12} \tilde g ^{ab} \tilde
 g^{cd} \tgr \tilde g_{cd,r} ,\nn\\
 \tilde b^r_{,r} &=& {\textstyle\frac12} \tilde g^{ac} \tilde g^{bd} \tilde
 g _{cd,r} \tilde g _{ab,r} - {\textstyle\frac12} \tilde g^{ab} \tgrr, \nn\\
 \tilde g_{ab} \gtrr &=& 2 \tilde g^{ac} \tilde g^{bd} \tgr \tilde
 g_{cd,r} - \tilde g^{ab} \tgrr , \nn\\
 \gtr \tilde g^{cd} _{,r} \left( \tilde g_{ab} \tilde g_{cd} +
 \tilde g_{ac} \tilde g_{bd} + \tilde g_{ad} \tilde g_{bc} \right)
 &=& \tilde g^{ab} \tilde g^{cd} \tgr \tilde g_{cd,r} + 2 \tilde
 g^{ab} \tilde g^{cd} \tilde g _{ac,r} \tilde g_{bd,r} .\nn
\end{eqnarray}
 Performing the integrations we obtain the contributions (modulo $\pi^{(m-1)/2}
 \sqrt { \tilde g}$)
 \begin{eqnarray*} \tilde \alpha_I &=&\textstyle \tilde g^{ab} \tilde g^{cd} \tgr \tilde
g_{cd,r}
 \left[ - {\textstyle\frac14} d_{2003} + {\textstyle\frac14} d_{0003} - \frac 3 2
 d_{2004} + {\textstyle\frac12} d_{0004} - 3 d_{2005} \right] \nn\\
 &+& \tilde g^{ab} \tilde g^{cd} \tilde g_{ac,r} \tilde g_{bd,r}
 \left[ d_{2003} - d_{0003} + 4 d_{2004} + d_{0004} - 6 d_{2005}
 \right] \nn\\
 &+& \tilde g^{ab} \tgrr \left[ - d_{2003} + {\textstyle\frac12} d_{0003} -
 2 d_{2004} \right] ,\nn\\
 \beta_I &=&\tilde g^{ab} \tilde g^{cd} \tgr \tilde g_{cd,r}
 \left[\textstyle\frac 3 {128} d_{0151} + {\textstyle\frac1{64}} d_{0151} + \frac{\sqrt{-1}} 8
 d_{1123} + \frac {\sqrt{-1}} {16} d_{1122} - \frac {\sqrt{-1}} 4 d_{1103}\right.\\&&\left.\qquad \textstyle- \frac {\sqrt{-1}} 8
 d_{1102} - {\textstyle\frac1{32}} d_{0111} \right] \nn\\
 &+&\textstyle\tilde g^{ab} \tilde g^{cd} \tilde g_{ac,r} \tilde g_{bd,r}
 \left[\frac 3 {64} d_{0151} + {\textstyle\frac1{32}} d_{0151} + \frac {\sqrt{-1}} 4
 d_{1123} - {\textstyle\frac18} d_{0131} + {\textstyle\frac18} d_{0111} \right] \nn\\
 &+&\tilde g^{ab} \tgrr \left[ {\textstyle\frac1{16}} d_{0131} - {\textstyle\frac18}
 d_{0111} \right] , \nn\\
 \gamma_I &=&\textstyle\tilde g^{ab} \tilde g^{cd} \tgr \tilde g_{cd,r}
 \left[ \frac 3 {128} d_{0241} + {\textstyle\frac1{64}} d_{0241} + \frac {\sqrt{-1}} 8
 d_{1213} + {\textstyle\frac1{32}} d_{0201} - \frac 3 {64} d_{0221}\right.\\&&\quad\left.\textstyle + \frac
 {\sqrt{-1}} {16} d_{1212} + {\textstyle\frac1{64}} d_{0221} \right] \nn\\
 &+&\textstyle\tilde g^{ab} \tilde g^{cd} \tilde g_{ac,r} \tilde g_{bd,r}
 \left[\frac 3 {64} d_{0241} + {\textstyle\frac1{32}} d_{0241} + \frac {\sqrt{-1}} 4
 d_{1213} - {\textstyle\frac18} d_{0221} + {\textstyle\frac18} d_{0201} \right] \nn\\
&+&\tilde g^{ab} \tgrr \left[{\textstyle\frac1{16}} d_{0221} - {\textstyle\frac18}
d_{0201} \right] \nn\\
\delta_I &=&\textstyle\tilde g^{ab} \tilde g^{cd} \tgr \tilde g_{cd,r}
 \left[\frac 5 {192} d_{0331} - {\textstyle\frac1{32}} d_{0311} \right]\nn\\
&+&\textstyle\tilde g^{ab} \tilde g^{cd} \tilde g_{ac,r} \tilde g_{bd,r}
 \left[\frac 5 {96} d_{0331} - {\textstyle\frac1{12}} d_{0311} \right] \nn\\
&+&\tilde g^{ab} \tgrr \left[{\textstyle\frac1{24}} d_{0311} \right] ,\nn\\
\epsilon_I &=&\tilde g^{ab} \tilde g^{cd} \tgr \tilde g_{cd,r}
 \left[{\textstyle\frac1{128}} d_{0421} \right] +
\tilde g^{ab} \tilde g^{cd} \tilde g_{ac,r} \tilde g_{bd,r}
 \left[{\textstyle\frac1{64}} d_{0421} \right].\nn
\end{eqnarray*}
Adding up all terms and simplifying using the functional equation
and the doubling formula for the $\Gamma$-function, the
contribution to the heat kernel coefficient reads
\begin{eqnarray}
& &(4\pi)^{-m/2} \textstyle\left\{\tilde g^{ab} \tilde g^{cd} \tgr \tilde
g_{cd,r} \frac{3 \alpha^2 - 16 \alpha - 27 }{384 (\alpha -6)}
\Gamma \left( \frac{1-\alpha} 2 \right) \right. \nn\\
& &\textstyle +\tilde g^{ab} \tilde g^{cd} \tilde g_{ac,r} \tilde g_{bd,r}
\frac{ 5 (9+4 \alpha - \alpha^2)}{192 (\alpha -6)} \Gamma \left(
\frac{ 1-\alpha} 2 \right)
\textstyle\left. +\tilde g^{ab} \tgrr \frac {\alpha +3 }{48} \Gamma
\left( \frac{ 1-\alpha } 2 \right)\right\}.\nn
\end{eqnarray}
This allows us to read off $A$, $B$, and $C$ from Equations (\ref{eqn-4.h}) and (\ref{eqn-4.i}) and to conclude:
\begin{eqnarray*} &&\kappa_\alpha ^3 =-
{\textstyle\frac1{12}}
(\alpha -1) \kappa_\alpha ,\quad
\kappa_\alpha^4 =\textstyle\frac{7-8\alpha + \alpha^2} {16 (\alpha -6)}
\kappa_\alpha ,\quad
\kappa_\alpha^5 =\frac{ 6 \alpha - 5 - \alpha^2}{8 (\alpha -6)
} \kappa_\alpha .\nn
\end{eqnarray*}
This completes the proof of Lemma \ref{lem-1.10}.\hfill\qedbox

\section{The proof of Theorem \ref{thm-1.6}}\label{sect-5}

\medbreak We adopt the notation of Theorem \ref{thm-1.6}.
 We may suppose that
the weighting function $F$ is supported in a boundary coordinate
neighborhood $U$ and that $F$ has the form $F=r^{-\alpha}G$ with
$G\in C_{\operatorname{comp}}^\infty(U)$. For appropriate functions
$q$ and $q^{bd}$ set
\begin{eqnarray*}
&&Op(q)f(y,r):=\int\int e^{\sqrt{-1}y\cdot
\omega+\sqrt{-1}r\tau}q(y,r,\omega,\tau,\lambda)\hat
f(\omega,\tau)\bar d\omega\bar d\tau,\\
&&Op^\prime(q^{bd})f(y,r):=\int\int e^{\sqrt{-1}y\cdot
\omega}q^{bd}(y,r,\omega,\tau,\lambda)\hat f(\omega,\tau)\bar
d\omega\bar d\tau,\\ &&\bar
d\omega:=(2\pi)^{1-m}d\omega,\quad\text{and}\quad\bar
d\tau:=(2\pi)^{-1}d\tau\,.
\end{eqnarray*}
Thus one has that $Op^\prime(q^{bd})=Op(e^{-\sqrt{-1}r\tau}q^{bd})$. In local coordinates
$(y,r)$, we follow the construction in \cite{Se69} to define the standard parametrix for $(\DD-\lambda)^{-1}$:
$$Q_N(\lambda)=\sum_{n=0}^N\left\{Op(q_{-2-n}(\lambda))-Op^\prime(q_{-2-n}^{bd}(\lambda))\right\}\,.$$
The $q_j$ and $q_j^{bd}$ are discussed in further detail in Section \ref{sect-4}; we refer in particular to Equations
(\ref{eqn-4.a}), (\ref{eqn-4.b}) and (\ref{eqn-4.c}). Note that the $q^{bd}_j$ here is the $h_j$ of Section \ref{sect-4}.

We pass to a parametrix $H_N$ for the operator $e^{-t\DD}$ by performing a
contour integration:
\begin{eqnarray*}
&&H_N(\DD,t):=\frac1{2\pi\sqrt{-1}}\int_\Gamma
e^{-t\lambda}Q_N(\lambda)d\lambda=\sum_{n=0}^N\{Op(e_{-2-n})-Op^\prime(e_{-2-n}^{bd})\}\,;
\end{eqnarray*}
in this expression, $e_j$ and $e_j^{bd}$ are the corresponding contour integrals of $q_j$ and $q_j^{bd}$. For another
construction of the heat parametrix see
\cite{Gr68}.

Take $\phi$ and $\psi$ with compact support in $U$, with $\psi\equiv1$ in a neighborhood of the support of $\phi$.
Then
$$R_N=\phi\left[e^{-t\DD}-H_N\right]\psi$$
has a kernel $k(y,r,y^\prime,r^\prime;t)$ which, for sufficiently large $N$, is $C^2$, $O(t^J)$ for large $J$
(depending on $N$), and vanishes when $r=0$, although not for $r^\prime=0$. To achieve vanishing in both variables, we
consider
\begin{equation}\label{eqn-5.a}
\begin{array}{l}
\phantom{=}\phi H_N(\DD, t/2)\psi^2H_N(\DD^*, t/2)^*\phi-\phi e^{-t\DD/2}\psi^2(e^{-t\DD^*/2})^*\phi\\
=-R_N\psi H_N(\DD^*, t/2)^*\phi-\phi H_N(\DD, t/2)\psi
R_N^*-R_NR_N^*\,.\vphantom{\vrule height 13pt}
\end{array}
\end{equation}
Here we compute the adjoints with respect to the measure $dydr$. The remainder given by Equation (\ref{eqn-5.a}) has a
kernel which is $C^2$, with derivatives $O(t^J)$ for large $J$, and vanishes when either $r=0$ or $r^\prime=0$. By
the pseudo-local properties of $H_N$ and $e^{-t\DD}$, the same is true of
$$\phi H_N(\DD, t/2)H_N(\DD^*, t/2)^*\phi-\phi e^{-t\DD}\phi\,.$$

Hence the kernel of this operator is $O(r^2t^J)$ for large $J$ and we can deduce an expansion for
$\Tr_{L^2}(Fe^{-t\DD})$ from an expansion of $\Tr_{L^2}(FH_NH_N^*\phi)$, taking $\phi\equiv1$ on the support of $F$.
The error will be holomorphic in $\alpha$ for $\operatorname{Re}(\alpha)<3$. Since $\Tr_{L^2}(Fe^{-t\DD})$ is also holomorphic
for
$\operatorname{Re}(\alpha)<3$, we may compute our expansion for $\operatorname{Re}(\alpha)<1$ and continue
analytically.

From the expansion of $q_j$ in powers of $(|\omega|^2+\tau^2-\lambda)^{-1}$ (See Section \ref{sect-4}) we get the estimates
\begin{eqnarray}\label{eqn-5.b}
&&|e_{-2-k}|\le Ct^{k/2}e^{-ct(|\omega|^2+\tau^2)}\\
&&\label{eqn-5.c}
|e_j^{bd}|\le C(|\omega|^2+\tau^2+1)^{-1}e^{-|\omega|^2t/2}
\end{eqnarray}
for suitably chosen constants $C$ and $c$. Set $\xi=(\omega,\tau)$. The $q_j$ and the $q_j^{bd}$ have
an appropriate homogeneity property \cite{Se69}. This homogeneity yields that
\begin{eqnarray}\label{eqn-5.d}
&&e_j(y,r,s\xi,t/s^2)=s^{2+j}e_j(y,r,\xi,t),\quad\text{and}\\
&&e_j^{bd}(y,r/s,s\xi,t/s^2)=s^{2+j}e_j^{bd}(y,r,\xi,t)\,.\label{eqn-5.e}
\end{eqnarray}
Furthermore, the kernel $k_1$ of $Op(e_j)Op(e_k^*)^*$ and the kernel $k_2$ of
$Op(e_j)Op^\prime(e_k^{bd*})^*$ on the diagonal are given, respectively, by:
\begin{eqnarray}\label{eqn-5.f}
&&k_1(y,r,y,r;t)=\int e_j(y,r,\xi,t/2)e_k(y,r,\xi,t/2)\bar d\xi\\
&&\qquad=(t/2)^{-(m+j+k+4)/2}\int e_j(y,r,\xi,1)e_k(y,r,\xi,1)\bar d\xi,\nonumber\\
\label{eqn-5.g} &&k_2(y,r,y,r;t)=(t/2)^{-(m+j+k+4)/2}\int
e_j(y,r,\xi,1)e_k^{bd}(y,r(2/t)^{1/2},\xi,1)\bar d\xi\,.
\end{eqnarray}
There are similar formulas for the kernel of $Op^\prime(e_j^{bd})Op(e_k^*)^*$ and for the kernel of
$Op^\prime(e_j^{bd})Op^\prime(e_k^{bd*})^*$ on the diagonal.

We integrate $F(y,r)k_i(y,r,y,r;t)$. An expansion of $Fe_je_k$ in
powers  $r^{\ell-\alpha}$ gives terms $t^{-(m+j+k+4)/2}$ times a
meromorphic function with simple poles at $\alpha=1,2,...$. For
other $\alpha$, this gives the interior terms in Theorem
\ref{thm-1.6}; the terms with $j+k$ odd vanish because of the parity
of $e_j$ and $e_k$ in $\xi$. For the terms in Equation
(\ref{eqn-5.g}), note that from Equations (\ref{eqn-5.c}) and
(\ref{eqn-5.e}), $e_k^{bd}$ decays exponentially as
$r\rightarrow\infty$ so that we may integrate in $r$ from $0$ to
$\infty$. An expansion of $Fe_j$ in powers of $r$ and a change of
variable $r/\sqrt{t}\rightarrow r$ gives boundary terms of the form
in Theorem \ref{thm-1.6} for $\alpha\ne1,2$.

This proves Theorem \ref{thm-1.6} for $\operatorname{Re}(\alpha)<1$ and the rest follows by analytic continuation.
Since for each $t$, the expansion is holomorphic in $\alpha$ for $\operatorname{Re}(\alpha)<3$, the residues arising from
the interior integrals at $\alpha=1,2$ must be cancelled by residues from the boundary terms. And, since the various
powers of $t$ are linearly independent, the residue of each interior coefficient must be cancelled by the residue
from the corresponding boundary coefficient. The residue from $t^n\Ireg(Fa_n)=t^n\Ireg(r^{-\alpha}Ga_n)$ is
\begin{equation}\label{eqn-5.h}
\left\{\begin{array}{lll}
-t^n{\mathcal{I}}^{bd}\{(Fa_n)_0\}&\text{at}&\alpha=1,\\
-t^n{\mathcal{I}}^{bd}\{(Fa_n)_1+(Fa_n)_0L_{aa}\}&\text{at}&\alpha=2\,.
\end{array}\right.\end{equation}
The corresponding boundary term when $\alpha=1$ is $t^{n-(\alpha-1)/2}a_{2n,\alpha}^{bd}(F,\DD)$. Let $R$ be the residue of
$a_{2n,\alpha}^{bd}(F,\DD)$ at $\alpha=1$. Then
$$t^{n-(\alpha-1)/2}a_{2n,\alpha}^{bd}(F,\DD)=\frac{t^nR}{\alpha-1}-\frac12t^n\ln(t)R+(const)t^n+O(\alpha-1)$$
Assertion (2) now follows. Since $t^nR$ must cancel the residue in
Equation (\ref{eqn-5.h}), we get the $\ln(t)$ coefficient in
Assertion (4) for $\alpha=1$. The proof of the rest of Assertion (4)
and also of Assertion (3) follows in a similar fashion.\hfill\qedbox

\section*{Acknowledgments} The authors acknowledge support by the Isaac Newton Institute in
Cambridge during the Programme on Spectral Theory and Partial
Differential Equations in July 2006 where this research began. The
research of M. van den Berg was supported by the London Mathematical
Society under Scheme 4, reference 4703, and by The Leverhulme Trust
, Research Fellowship 2008/0368. The research of P. Gilkey was
supported by the Max Planck Institute for Mathematics in the
Sciences (Germany) and by Project MTM2006-01432 (Spain). The
research of K. Kirsten was supported by the Baylor University Summer
Sabbatical Program and by the University Research Council.

\end{document}